\documentclass[
a4paper,
11pt
]{article}

\usepackage{
amssymb,
amsmath,
amsthm,
eucal,
}
\makeatletter
 
  \@addtoreset{equation}{section}
 \makeatother
\theoremstyle{plain}

\newtheorem{theorem}{Theorem}[section]
\newtheorem{lemma}[theorem]{Lemma}
\newtheorem{proposition}[theorem]{Proposition}
\newtheorem{corollary}[theorem]{Corollary}
\theoremstyle{definition}
\newtheorem{definition}[theorem]{Definition}
\newtheorem{remark}[theorem]{Remark}

\newtheorem{assumption}[theorem]{Assumption}
\newtheorem{acknowledgements}{Acknowledgements\!\!}

\newcommand{\bignorm}[1]{{\left|\left|#1\right|\right|}}
\newcommand{\norm}[1]{{||#1||}}
\newcommand{\abs}[1]{{\left|#1\right|}}

\newcommand{\dderiv}[4]{{\partial_{#1}^{#3}\partial_{#2}^{#4}}}

\newcommand{\what}[1]{{\widehat{#1}}}\newcommand{\wtilde}[1]{{\widetilde{#1}}}

\def\supp{\mathop{\mathrm{supp}}\nolimits}

\def\Id{\mathop{\mathrm{Id}}\nolimits}

\def\Op{\mathop{\mathrm{Op}}\nolimits}
\def\IK{\mathop{\mathrm{IK}}\nolimits}
\def\WKB{\mathop{\mathrm{WKB}}\nolimits}

\def\Im{\mathop{\mathrm{Im}}\nolimits}
\def\R{{\mathbb{R}}}
\def\Z{{\mathbb{Z}}}
\def\N{{\mathbb{N}}}
\def\C{{\mathbb{C}}}
\def\Sphere{{\mathbb{S}}}
\def\S{{\mathcal{S}}}

\def\B{{\mathcal{B}}}

\def\X{{\mathcal{X}}}
\def\Y{{\mathcal{Y}}}

\def\<{{\langle}}
\def\>{{\rangle}}
\def\ep{{\varepsilon}}

\title%[Strichartz estimates with unbounded potentials]
{Strichartz estimates for Schr\"odinger equations with variable coefficients and unbounded potentials}
\author%[Haruya Mizutani]
{Haruya Mizutani${}^*$}
\date{\empty}
\begin{document}
\maketitle

\footnotetext{ %2000 MSC numbers
}
\footnotetext{ %key words and phrases
}
\footnotetext{${}^*$Research Institute for Mathematical Sciences, Kyoto University, Kyoto 606-8502, Japan. Current address: Department of Mathematics, Gakushuin University, 1-5-1 Mejiro, Toshima-ku, Tokyo 171-8588, Japan. E-mail: \texttt{haruya@math.gakushuin.ac.jp}. The author was partially supported by GCOE `Fostering top leaders in mathematics', Kyoto University. 
}

\begin{abstract}
The present paper is concerned with Schr\"odinger equations with variable coefficients and unbounded electromagnetic potentials, where the kinetic energy part is a long-range perturbation of the flat Laplacian and the electric (resp. magnetic) potential can grow subquadratically (resp. sublinearly) at spatial infinity. 
We prove sharp (local-in-time) Strichartz estimates, outside a large compact ball centered at origin, for any admissible pair including the endpoint. 
Under the nontrapping condition on the Hamilton flow generated by the kinetic energy, global-in-space estimates are also studied. 
Finally, under the nontrapping condition, we prove Strichartz estimates with an arbitrarily small derivative loss without asymptotic flatness on the coefficients. 
\end{abstract}

%%%%%%%%%%%%%%%%%%%%%%%%%%%%%%%%%%%%%%%%%%%%%%%%%%%%%%%%%%%%%			%Introduction			%%%%%%%%%%%%%%%%%%%%%%%%%%%%%%%%%%%%%%%%%%%%%%%%%%%%%%%%%%%%%%%%%%%%%
\section{Introduction}
\label{section_introduction}
In this paper, we study sharp (local-in-time) Strichartz estimates for Schr\"odinger equations with variable coefficients and unbounded electromagnetic potentials. 
More precisely, we consider the following Schr\"odinger operator:
$$
H=\frac12\sum_{j,k=1}^d (-i\partial_j-A_j(x)) g^{jk}(x) (-i\partial_k-A_k(x))+V(x),\quad x\in\R^d,
$$
where $d\ge1$ is the spatial dimension. 
Throughout the paper we assume that $g^{jk},V$ and $A_j$ are smooth real-valued functions on $\R^d$ and that $(g^{jk}(x))_{j,k}$ is symmetric and positive definite:
$$
\sum_{j,k=1}^d g^{jk}(x)\xi_j\xi_k \ge c|\xi|^2,\quad x,\xi\in\R^d,
$$
with some $c>0$. 
Moreover, we suppose the following condition:

\begin{assumption}
\label{assumption_A}
There exists $\mu\ge0$ such that for any $\alpha\in \Z^d_+$,
\begin{align*}
|\partial_x^\alpha(g^{jk}(x)-\delta_{jk})|&\le C_{\alpha}\<x\>^{-\mu-|\alpha|},\\
|\partial_x^\alpha A_j(x)|&\le C_{\alpha}\<x\>^{1-\mu-|\alpha|},\\
|\partial_x^\alpha V(x)|&\le C_{\alpha}\<x\>^{2-\mu-|\alpha|},\quad x\in\R^d.
\end{align*}
\end{assumption}

Then, it is well known that $H$ admits a unique self-adjoint realization on $L^2(\R^d)$, which we denote by the same symbol $H$. By the Stone theorem, $H$ generates a unique unitary propagator $e^{-itH}$ on $L^2(\R^d)$ such that the solution to the Schr\"odinger equation:
$$
\label{Schrodinger_equation}
i\partial_t u(t)=Hu(t),\ t\in\R;\quad u|_{t=0}=\varphi \in L^2(\R^d),
$$
is given by $u(t)=e^{-itH}\varphi$.

In order to explain the purpose of the paper, we recall some known results. Let us first recall well known properties of the free propagator $e^{-itH_0}$, where $H_0=-\Delta/2$. The distribution kernel of $e^{-itH_0}$ is given explicitly by $(2\pi it)^{-d/2}e^{i|x-y|^2/(2t)}$ and $e^{-itH_0}\varphi$ thus satisfies the dispersive estimate:
$$
\norm{e^{-itH_0}\varphi}_{L^\infty(\R^d)} \le C|t|^{-d/2}\norm{\varphi}_{L^1(\R^d)},\quad t\neq0. 
$$
Moreover, $e^{-itH_0}$ enjoys the following (global-in-time) Strichartz estimates:
$$
\norm{e^{-itH_0}\varphi}_{L^p(\R;L^q(\R^d))} \le C\norm{\varphi}_{L^2(\R^d)},
$$
where $(p,q)$ satisfies the following admissible condition:
\begin{align}
\label{admissible}%\label{admissible}
p\ge2,\quad \frac2p=d\left(\frac12-\frac1q\right),\quad  (d,p,q) \neq (2,2,\infty).
\end{align}
Strichartz estimates imply that, for any $\varphi\in L^2$, $e^{-itH_0}\varphi\in \bigcap_{q\in Q_d} L^q$ for a.e. $t\in\R$, where $Q_1=[2,\infty]$, $Q_2=[2,\infty)$ and $Q_d=[2,2d/(d-2)]$ for $d\ge3$. These estimates hence can be regarded as $L^p$-type smoothing properties of Schr\"odinger equations, and have been widely used in the study of nonlinear Schr\"odinger equations (see, \emph{e.g.}, \cite{Cazenave}). Strichartz estimates for $e^{-itH_0}$ were first proved by Strichartz \cite{Strichartz} for a restricted pair of $(p,q)$ with $p=q=2(d+2)/d$, and have been generalized for $(p,q)$ satisfying \eqref{admissible} and $p\neq2$ by \cite{Ginibre_Velo}. The endpoint estimate $(p,q)=(2,2d/(d-2))$ for $d\ge3$ was obtained by \cite{Keel_Tao}. 

For Schr\"odinger operators with electromagnetic potentials, \emph{i.e.}, $H=\frac12(-i\partial_x-A)^2+V$, (short-time) dispersive and (local-in-time) Strichartz estimates have been extended with potentials decaying at infinity \cite{Yajima1} or growing at infinity \cite{Fujiwara,Yajima2}. In particular, it was shown by  \cite{Fujiwara,Yajima2} that if $g^{jk}=\delta_{jk}$, $V$ and $A$ satisfy Assumption \ref{assumption_A} with $\mu\ge0$ and all derivatives of the magnetic field $B=dA$ are of short-range type, then $e^{-itH}\varphi$ satisfies (short-time) dispersive estimates:
$$
\norm{e^{-itH}\varphi}_{L^\infty(\R^d)} \le C|t|^{-d/2}\norm{\varphi}_{L^1(\R^d)},
$$
for sufficiently small $t\neq0$. Local-in-time Strichartz estimates, which have the forms
$$
\norm{e^{-itH}\varphi}_{L^p([-T,T];L^q(\R^d))} \le C_T\norm{\varphi}_{L^2(\R^d)},\quad T>0,
$$ 
are immediate consequences of this estimate and the $TT^*$-argument due to Ginibre-Velo \cite{Ginibre_Velo} (see Keel-Tao \cite{Keel_Tao} for the endpoint estimate). For the case with singular electric potentials or with supercritical electromagnetic potentials, we refer to \cite{Yajima1,Yajima3,Yajima_Zhang,Dancona_Fanelli} and reference therein. We mention that global-in-time dispersive and Strichartz estimates for scattering states have been also studied under suitable decaying conditions on potentials and assumptions for zero energy; see \cite{JSS,Yajima5,Schlag1,EGS,D'AFVV} and reference therein. We also mention that there is no result on sharp global-in-time dispersive estimates for magnetic Schr\"odinger equations.

On the other hand, the influence of the geometry on the behavior of solutions to linear and nonlinear partial differential equations has been extensively studied. From this geometric viewpoint, sharp local-in-time Strichartz estimates for Schr\"odinger equations with variable coefficients (or, more generally, on manifolds) have recently been investigated by many authors under several conditions on the geometry; see, \emph{e.g.}, \cite{Staffilani_Tataru,BGT,Robbiano_Zuily,HTW,Bouclet_Tzvetkov_1,Bouclet,BGH,Mizutani1} and reference therein. In \cite{Staffilani_Tataru}, \cite{Robbiano_Zuily}, \cite{Bouclet_Tzvetkov_1}, the authors studied the case on the Euclidean space with nontrapping asymptotically flat metrics. 
The case on the nontrapping asymptotically conic manifold was studied by \cite{HTW} and \cite{Mizutani1}. 
In \cite{Bouclet} the author considered the case of nontrapping asymptotically hyperbolic manifold. 
For the trapping case, it was shown in \cite{BGT} that Strichartz estimates with a loss of derivative $1/p$ hold on any compact manifolds without boundaries. 
They also proved that the loss $1/p$ is optimal in the case of $M=\Sphere^d$. In \cite{Bouclet_Tzvetkov_1}, \cite{Bouclet} and \cite{Mizutani1}, the authors proved sharp Strichartz estimates, outside a large compact set, without the nontrapping condition. 
More recently, it was shown in \cite{BGH} that sharp Strichartz estimates still hold for the case with hyperbolic trapped trajectories of sufficiently small fractal dimension. 
We mention that there are also several works on global-in-time Strichartz estimates in the case of long-range perturbations of the flat Laplacian on $\R^d$ (\cite{Bouclet_Tzvetkov_2,Tataru,MMT}).

While (local-in-time) Strichartz estimates are well studied subjects for both of these two cases (at least under the nontrapping condition), the literature is more sparse for the mixed case. In this paper we give a unified approach to a combination of these two kinds of results. More precisely, under Assumption \ref{assumption_A} with $\mu>0$, 
we prove (1) sharp local-in-time Strichartz estimates, outside a large compact set centered at origin, without the nontrapping condition; (2) the global-in-space estimates with the nontrapping condition. Under the nontrapping condition and Assumption \ref{assumption_A} with $\mu\ge0$, we also show local-in-time Strichartz estimates with an arbitrarily small derivative loss. We mention that all results include the endpoint estimates $(p,q)=(2,2d/(d-2))$ for $d\ge3$. This is a natural continuation of author's previous work \cite{Mizutani2}, which was concerned with the non-endpoint estimates for the case with at most linearly growing potentials.

In the sequel, $F(*)$ denotes the characteristic function designated by $(*)$. We now state the main result. 

%theorem
\begin{theorem}[Strichartz estimates near infinity]
\label{theorem_1}
Suppose that $H$ satisfies Assumption \ref{assumption_A} with $\mu>0$. Then, there exists $R_0>0$ such that for any $T>0$, $p\ge2$, $q<\infty$, $2/p=d(1/2-1/q)$ and $R\ge R_0$, we have
\begin{align}
\label{theorem_1_1}
\norm{F(|x|>R)e^{-itH}\varphi}_{L^p([-T,T];L^q(\R^d))} \le C_T \norm{\varphi}_{L^2(\R^d)},
\end{align}
where $C_T>0$ may be taken uniformly with respect to $R$.
\end{theorem}

To state the result on global-in-space estimates, we recall the \emph{nontrapping condition}. Let us denote by $k(x,\xi)$ the classical kinetic energy:
$$
k(x,\xi)=\frac12\sum_{j,k=1}^dg^{jk}(x)\xi_j\xi_k,
$$
and by $(y_0(t,x,\xi),\eta_0(t,x,\xi))$ the Hamilton flow generated by $k(x,\xi)$:
$$
\dot{y}_0(t)=\partial_\xi k(y_0(t),\eta_0(t)),\ \dot{\eta}_0(t)=-\partial_x k(y_0(t),\eta_0(t));\quad
(y_0(0),\eta_0(0))=(x,\xi). 
$$
Note that the Hamiltonian vector field $H_{k}=\partial_\xi k\cdot\partial_x-\partial_x k\cdot\partial_\xi$ generated by $k$ is complete on $\R^{2d}$ since $(g^{jk})$ satisfies the uniform elliptic condition. Hence, $(y_0(t,x,\xi),\eta_0(t,x,\xi))$  exists for all $t \in \R$. 

%{definition}
\begin{definition}
We say that $k(x,\xi)$ satisfies the nontrapping condition if for any $(x,\xi) \in \R^{2d}$ with $\xi \neq 0$,
\begin{align}
\label{nontrapping}
|y_0(t,x,\xi)| \to +\infty\ as\ t \to \pm \infty. 
\end{align}
\end{definition}
To control the asymptotic behavior of the flow, we also impose the following condition which is the classical analogue of Mourre's inequality:

%assumption
\begin{assumption}[Convexity near infinity]
\label{assumption_B}
There exists $f\in C^\infty(\R^d)$ satisfying $f\ge1$ and $f\to +\infty$ as $|x|\to+\infty$ such that $\partial^\alpha f\in L^\infty(\R^d)$ for any $|\alpha|\ge2$ and 
$$
H_k (H_kf)(x,\xi)\ge ck(x,\xi)
$$
on $\{(x,\xi)\in\R^{2d};\ f(x)\ge R\}$ for some positive constants $c,R>0$. \end{assumption}

Note that if $|\partial_x g^{jk}(x)|=o(|x|^{-1})$ as $|x|\to+\infty$, then Assumption \ref{assumption_B} holds with $f(x)=1+|x|^2$. In particular, Assumption \ref{assumption_A} with $\mu>0$ implies Assumption \ref{assumption_B}. Moreover, if $g^{jk}(x)=(1+a_1\sin(a_2\log r))\delta_{jk}$ for $a_1\in\R,a_2>0$ with $a_1^2(1+a_2^2)<1$ and for $r=|x|\gg1$, then Assumption \ref{assumption_B} holds with $f(r)=(\int_0^r(1+a_1\sin(a_2\log t))^{-1}dt)^2$. For more examples, we refer to \cite[Section 2]{Doi}. 

The second result then is the following.

%theorem
\begin{theorem}	[Global-in-space Strichartz estimates]		
\label{theorem_2}
Suppose that $H$ satisfies Assumption \ref{assumption_A} with $\mu\ge0$. Let $T>0$, $p\ge2$, $q<\infty$ and $2/p=d(1/2-1/q)$. Then, for any $r>0$, there exists $C_{T,r}>0$ such that
\begin{align}
\label{theorem_2_1}
\norm{F(|x|<r)e^{-itH}\varphi}_{L^p([-T,T];L^q(\R^d))} \le C_{T,r} \norm{\<H\>^{\frac{1}{2p}}\varphi}_{L^2(\R^d)}.
\end{align}
Moreover if we assume in addition that $k(x,\xi)$ satisfies the nontrapping condition \eqref{nontrapping} and that Assumption \ref{assumption_B}, then 
\begin{align}
\label{theorem_2_2}
\norm{F(|x|<r)e^{-itH}\varphi}_{L^p([-T,T];L^q(\R^d))} \le C_{T,r} \norm{\varphi}_{L^2(\R^d)}.
\end{align}
 In particular, combining with Theorem \ref{theorem_1}, we have (global-in-space) Strichartz estimates
$$
\norm{e^{-itH}\varphi}_{L^p([-T,T];L^q(\R^d))} \le C_{T} \norm{\varphi}_{L^2(\R^d)},
$$
under the nontrapping condition \eqref{nontrapping}, provided that $\mu>0$. 
\end{theorem}

For the general case we have the following  partial result.

%theorem
\begin{theorem}[Near sharp estimates without asymptotic flatness]		
\label{theorem_3}Suppose that $H$ satisfies Assumption \ref{assumption_A} with $\mu\ge0$ and $k(x,\xi)$ satisfies the nontrapping condition \eqref{nontrapping}. Assume also Assumption \ref{assumption_B}. Let $T>0$, $p\ge2$, $q<\infty$ and $2/p=d(1/2-1/q)$. Then, for any $\ep>0$, there exists $C_{T,\ep}>0$ such that
$$
\norm{e^{-itH}\varphi}_{L^p([-T,T];L^q(\R^d))} \le C_{T,\ep} \norm{\<H\>^{\ep}\varphi}_{L^2(\R^d)}. 
$$
\end{theorem}

There are some remarks.

%remark
\begin{remark}
(1) The estimates of forms \eqref{theorem_1_1}, \eqref{theorem_2_1} and \eqref{theorem_2_2} have been proved by \cite{Staffilani_Tataru,Bouclet_Tzvetkov_1} when $A\equiv0$ and $V$ is of long-range type. Theorems \ref{theorem_1} and \ref{theorem_2} hence are regarded as generalizations of their results for the case with growing electromagnetic potential perturbations.\\
(2) The only restriction for admissible pairs, in comparison to the flat case, is to exclude  $(p,q)=(4,\infty)$ for $d=1$, which is due to the use of the Littlewood-Paley decomposition.\\
(3) The missing derivative loss $\<H\>^{\ep}$ in Theorem \ref{theorem_3} is due to the use of the following local smoothing effect (due to Doi \cite{Doi}):
$$
\norm{\<x\>^{-1/2-\ep}\<D\>^{1/2}e^{-itH}\varphi}_{L^2([-T,T];L^2(\R^d))} \le C_{T,\ep}\norm{\varphi}_{L^2(\R^d)}.
$$
It is well known that this estimate does not holds when $\ep=0$ even for $H=H_0$. We would expect that Theorem \ref{theorem_1} still holds true for the case with critical electromagnetic potentials in the following sense:
$$
\<x\>^{-1}|\partial_x^\alpha A_j(x)|+
\<x\>^{-2}|\partial_x^\alpha V(x)| \le C_{\alpha\beta}\<x\>^{-|\alpha|}, 
$$
(at least if $g^{jk}$ satisfies the bounds in Assumption \ref{assumption_A} with $\mu>0$). However, this is beyond our techniques (see also remark \ref{remark_IK}).
\end{remark}

%%%%%%%%%%%%%%%%%%%%%%%%%%%%%%%%%%%%%%%%%%%%%%%			outline of the paper							       %%%%%%%%%%%%%%%%%%%%%%%%%%%%%%%%%%%%%%%%%%%%%%%%

The rest of the paper is devoted to the proofs of Theorems \ref{theorem_1}, \ref{theorem_2} and \ref{theorem_3}. Throughout the paper we use the following notations: $\<x\>$ stands for $\sqrt{1+|x|^2}$. We write $L^q=L^q(\R^d)$ if there is no confusion. For Banach spaces $X$ and $Y$, we denote by $\norm{\cdot}_{X\to Y}$ the operator norm from $X$ to $Y$. We write $\Z_+=\N\cup\{0\}$ and denote the set of multi-indices by $\Z^d_+$. We denote by $K$ the kinetic energy part of $H$ and by $H_0$ the free Schr\"odinger operator:
$$
K=-\frac12\sum_{j,k=1}^d\partial_jg^{jk}(x)\partial_k,\quad H_0=-\frac12\Delta=-\frac12\sum_{j=1}^d\partial_j^2.
$$ 
We set two symbols $p(x,\xi)$ and $p_1(x,\xi)$ defined by
\begin{equation}
\begin{aligned}
\label{symbols_1}
p(x,\xi)=\frac12\sum_{j,k=1}^dg^{jk}(x)(\xi_j-A_j(x))(\xi_k-A_k(x))+V(x),\\
p_1(x,\xi)=-\frac i2\sum_{j,k=1}^d\left(\frac{\partial g^{jk}}{\partial x_j}(x)(\xi_k-A_k(x))-g^{jk}(x)\frac{\partial A_k}{\partial x_j}(x)\right).
\end{aligned}
\end{equation}
Note that Assumption \ref{assumption_A} implies
\begin{equation}
\begin{aligned}
\label{estimates_on_symbols_0}
|\dderiv{x}{\xi}{\alpha}{\beta}p(x,\xi)|
&\le C_{\alpha\beta}\<x\>^{-|\alpha|}\<\xi\>^{-|\beta|}(|\xi|^2+\<x\>^{2-\mu}),\\
|\dderiv{x}{\xi}{\alpha}{\beta}p_1(x,\xi)|
&\le C_{\alpha\beta}\<x\>^{-|\alpha|}\<\xi\>^{-|\beta|}(\<x\>^{-1-\mu}|\xi|+\<x\>^{-\mu}).
\end{aligned}
\end{equation} 
For $h\in (0,1]$ we consider $H^h:=h^2H$ as a semiclassical Schr\"odinger operator with $h$-dependent electromagnetic potentials $h^2V$ and $hA_j$. The corresponding symbols $p_h$ and $p_{1,h}$ are also defined by
\begin{equation}
\begin{aligned}
\label{symbols_2}
p_h(x,\xi)&=\frac12\sum_{j,k=1}^dg^{jk}(x)(\xi_j-hA_j(x))(\xi_k-hA_k(x))+h^2V(x),\\
p_{1,h}(x,\xi)&=-\frac i2\sum_{j,k=1}^d\left(\frac{\partial g^{jk}}{\partial x_j}(x)(\xi_k-hA_k(x))-hg^{jk}(x)\frac{\partial A_k}{\partial x_j}(x)\right).
\end{aligned}
\end{equation}
It is easy to see that $H=\Op(p)+\Op(p_1)$ and $H^h=\Op_h(p_h)+h\Op_h(p_{1,h})$. 

Before starting the details of the proofs, we here describe the main ideas. At first we note that, since our Hamiltonian $H$ is not bounded below, the Littlewood-Paley decomposition associated with $H$ seems to be false for $p\neq2$ in general. To overcome this difficulty, we consider the following partition of unity on the phase space $\R^{2d}$:
$$
\psi_\ep(x,\xi)+\chi_\ep(x,\xi)=1,
$$
where $\psi_\ep$ is supported in $\{(x,\xi);\<x\><\ep|\xi|\}$ for some sufficiently small constant $\ep>0$. It is easy to see that the symbol $p(x,\xi)$ is elliptic on $\supp \psi_\ep$:
$$
C^{-1}|\xi|^2 \le p(x,\xi) \le C|\xi|^2,\quad (x,\xi)\in\supp\psi_\ep,
$$
 and we hence can prove a Littlewood-Paley type decomposition of the following form:
\begin{align*}
\norm{\Op(\psi_\ep)u}_{L^q}
\le C_q \norm{u}_{L^2}
+C_q\bigg(\sum_{h=2^{-j}, j\ge0}\norm{\Op_h(a_h)f(h^2 H)u}_{L^q}^2\bigg)^{1/2}, 
\end{align*}
where $2\le q<\infty$, $\{f(h^2\cdot);h=2^{-j},j\ge0\}$ is a $4$-adic partition of unity on $[1,\infty)$ and $a_h$ is an appropriate $h$-dependent symbol supported in $\{|x|<1/h,\ |\xi|\in I\}$ for some open interval $I\Subset(0,\infty)$, $\Op(\psi_\ep)$ and $\Op_h(a_h)$ denote the corresponding pseudodifferential and semiclassical pseudodifferential operators, respectively. 

Then, the idea of the proof of Theorem \ref{theorem_1} is as follows. In view of the above Littlewood-Paley estimate, the proof is reduced to that of Strichartz estimates for $F(|x|>R)\Op_h(a_h)e^{-itH}$ and $\Op(\chi_\ep)e^{-itH}$. In order to prove Strichartz estimates for $F(|x|>R)\Op_h(a_h)e^{-itH}$, we use semiclassical approximations of Isozaki-Kitada type. We however note that because of the unboundedness of potentials with respect to $x$, it is difficult to construct directly such approximations. To overcome this difficulty, we introduce a modified Hamiltonian $\wtilde{H}$ due to \cite{Yajima_Zhang} so that $\wtilde{H}=H$ for $|x|\le L/h$ and $\wtilde{H}=K$ for $|x|\ge 2L/h$ for some constant $L\ge1$. Then, $\wtilde{H}^h=h^2\wtilde{H}$ can be regarded as a ``long-range perturbation" of the semiclassical free Schr\"odinger operator $H_0^h=h^2H_0$. We also introduce the corresponding modified symbol $\wtilde{p}_h(x,\xi)$ so that $\wtilde{p}_h(x,\xi)=p_h(x,\xi)$ for $|x|\le L/h$ and $\wtilde{p}_h(x,\xi)=k(x,\xi)$ for $|x|\ge 2L/h$. Let $a^\pm_h$ be supported in outgoing and incoming regions
$
\{R<|x|<1/h,\ |\xi|\in I,\ \pm\hat{x}\cdot\hat{\xi}>1/2\}
$, 
respectively, so that $F(|x|>R)a_h=a_h^++a_h^-$, where $\hat{x}=x/|x|$. Rescaling $t\mapsto th$, we first construct  the semiclassical approximations for $e^{-it\wtilde{H}^h/h}\Op_h(a^\pm_h)^*$ of the following forms
$$
e^{-it\wtilde{H}^h/h}\Op_h(a_h^\pm)^*=J_h(S^\pm_h,b^\pm_h)e^{-itH_0^h/h}J_h(S^\pm_h,c^\pm_h)^*+O(h^N),\quad 0\le \pm t\le 1/h,
$$
respectively, where $S^\pm_h$ solve the Eikonal equation associated to $\wtilde{p}_h$ and $J_h(S^\pm_h,b^\pm_h)$ and $J_h(S^\pm_h,c^\pm_h)$ are associated semiclassical Fourier integral operators. The method of the construction is similar to as that of Robert \cite{Robert}. On the other hand, we will see that if $L\ge1$ is large enough, then the Hamilton flow generated by $\wtilde{p}_h$ with initial conditions in $\supp a_h^\pm$ cannot escape from $\{|x|\le L/h\}$ for $0<\pm t\le 1/h$, respectively, \emph{i.e.}, 
$$
\pi_x\left(\exp tH_{\wtilde{p}_h}(\supp a_h^\pm)\right)\subset \{|x|\le L/h\},\quad 0<\pm t\le 1/h.
$$
Since $\wtilde{p}_h=p_h$ for $|x|\le L/h$, we have 
$$
\exp tH_{\wtilde{p}_h}(\supp a_h^\pm)=\exp tH_{p_h}(\supp a_h^\pm),
\quad0<\pm t\le 1/h.
$$ 
We thus can expect (at least formally) that the corresponding two quantum evolutions  are approximately equivalent modulo some smoothing operator. We will prove the following rigorous justification of this formal consideration:
$$
\norm{(e^{-itH^h/h}-e^{-it\wtilde{H}^h/h})\Op_h(a_h^\pm)^*}_{L^2 \to L^2} \le C_{M} h^{M},\quad 0\le\pm t\le 1/h,\ M\ge0,
$$
where $H^h=h^2H$. By using such approximations for $e^{-itH^h/h}\Op_h(a^\pm_h)^*$, we prove local-in-time dispersive estimates for $\Op_h(a^\pm_h)e^{-itH}\Op_h(a^\pm_h)^*$:
$$
\norm{\Op_h(a^\pm_h)e^{-itH}\Op_h(a^\pm_h)^*}_{L^1 \to L^\infty} \le C |t|^{-d/2},\quad 0<h\ll1,\ 0<|t|<1.
$$
 Strichartz estimates follow from these estimates and the abstract Theorem due to Keel-Tao \cite{Keel_Tao}.
 
Strichartz estimates for $\Op(\chi_\ep)e^{-itH}$ follow from the following short-time dispersive estimate:
$$
\norm{\Op(\chi_\ep)e^{-itH}\Op(\chi_\ep)^*}_{L^1 \to L^\infty} \le C_\ep |t|^{-d/2},\quad 0<|t|<t_\ep\ll1.
$$ 
To prove this, we first construct an approximation for $e^{-itH}\Op(\chi_\ep)^*$ of the following form:
$$
e^{-itH}\Op(\chi_\ep)^*=J(\Psi,a)+O_{H^{-\gamma}\to H^\gamma}(1),\quad |t|<t_\ep,\ \gamma>d/2,
$$
where the phase function $\Psi=\Psi(t,x,\xi)$ is a solution to a time-dependent Hamilton-Jacobi equation associated to $p(x,\xi)$ and $J(\Psi,a)$ is the corresponding Fourier integral operator. In the construction, the following fact plays an important rule:
$$
|\dderiv{x}{\xi}{\alpha}{\beta}p(x,\xi)|\le C_{\alpha\beta},\quad (x,\xi)\in \supp \chi_\ep,\ |\alpha+\beta|\ge2.
$$
We note that if $(g^{jk})_{jk}-\Id_d\neq0$ depends on $x$ then these bounds do not hold without such a restriction of the support. Using these bounds, we construct the phase function $\Psi(t,x,\xi)$ such that
$$
|\dderiv{x}{\xi}{\alpha}{\beta}(\Psi(t,x,\xi)-x\cdot\xi+p(x,\xi))|
\le C_{\alpha\beta}|t|^2\<x\>^{2-|\alpha+\beta|}.
$$
Then, we can follow a classical argument (due to \cite{Kitada_Kumanogo}) and construct the FIO $J(\Psi,a)$. By the composition formula, $\Op(\chi_\ep)J(\Psi,a)$ is also a FIO and dispersive estimates for this operator follow from the standard stationary phase method. Finally, using an Egorov type lemma, we prove that the remainder, $\Op(\chi_\ep)(e^{-itH}\Op(\chi_\ep)^*-J(\Psi,a))$, has a smooth kernel for sufficiently small $t$.

The  proof of Theorem \ref{theorem_2} is based on a standard idea by \cite{Staffilani_Tataru}, see also \cite{BGT,Bouclet_Tzvetkov_1}. Strichartz estimates with loss of derivatives $\<H\>^{1/(2p)}$ follow from semiclassical Strichartz estimates up to time scales of order $h$, which can be verified by the standard argument. Moreover, under the nontrapping condition, we will prove that the missing $1/p$ derivative loss can be recovered by using local smoothing effects due to Doi \cite{Doi}.

The  proof of Theorem \ref{theorem_3} is based on a slight modification of that of Theorem \ref{theorem_2}. By virtue of the Strichartz estimates for $\Op(\chi_\ep)e^{-itH}$ and the Littlewood-Paley decomposition, it suffices to show 
$$
\norm{\Op_h(a_h)e^{-itH}\varphi}_{L^p([-T,T];L^q)} \le C_Th^{-\ep}\norm{\varphi}_{L^2}, \quad 0<h\ll1.
$$
To prove this estimate, we first prove semiclassical Strichartz estimates for $e^{-itH}\Op_h(a_h)^*$ up to time scales of order $hR$, where $R=\inf |\pi_x(\supp a_h)|$. The proof is based on a refinement of the standard WKB approximation for the semiclassical propagator $e^{-itH^h/h}\Op_h(a_h)^*$. Combining semiclassical Strichartz estimates with a partition of unity argument with respect to $x$, we will obtain the following Strichartz estimate with an inhomogeneous error term:
$$
\norm{\Op_h(a_h)e^{-itH}\varphi}_{L^p([-T,T];L^q)}
\le C_T\norm{\varphi}_{L^2}+C\norm{\<x\>^{-1/2-\ep}h^{-1/2-\ep}\Op_h(a_h)e^{-itH}\varphi}_{L^2([-T,T];L^2)},
$$
for any $\ep>0$, which, combined with local smoothing effects, implies Theorem \ref{theorem_3}. 

The paper is organized as follows. We first record some known results on the semiclassical pseudodifferential calculus and prove the above Littlewood-Paley decomposition in Section \ref{section_functional_calculus}. Using dispersive estimates, which will be studied in Sections \ref{section_IK} and \ref{section_WKB}, we shall prove  Theorem \ref{theorem_1} in Section \ref{section_proof_theorem_1}. We construct approximations of Isozaki-Kitada type and prove dispersive estimates for $\Op_h(a^\pm_h)e^{-itH}\Op_h(a^\pm_h)^*$ in Section \ref{section_IK}. Section \ref{section_WKB}  discuss the dispersive estimates for $\Op(\chi_\ep)e^{-itH}\Op(\chi_\ep)^*$. The proof of Theorem \ref{theorem_2} and Theorem \ref{theorem_3} are given in Section \ref{section_proof_theorem_2} and Section \ref{section_proof_theorem_3}, respectively.

\begin{acknowledgements}
The author would like to express his sincere thanks to Professor Erik Skibsted for valuable discussions and for hospitality at Institut for Matematiske Fag, Aarhus Universitet, where a part of this work was carried out. He also would like to thank the referee for very careful reading the manuscript and for providing valuable suggestions, which substantially helped improving the quality of the paper. 
\end{acknowledgements}

%%%%%%%%%%%%%%%%%%%%%%%%%%%%%%%%%%%%%%%%%%%%%%%%%%%%%%			%Preliminaries			%%%%%%%%%%%%%%%%%%%%%%%%%%%%%%%%%%%%%%%%%%%%%%%%%%%%%%%

\section{Semiclassical functional calculus}
\label{section_functional_calculus}
Throughout this section we assume Assumption \ref{assumption_A} with $\mu\ge0$, \emph{i.e.},
\begin{align}
\label{assumption_A'}
|\partial_x^\alpha g^{jk}(x)|
+\<x\>^{-1}|\partial_x^\alpha A_j(x)|+
\<x\>^{-2}|\partial_x^\alpha V(x)| \le C_{\alpha\beta}\<x\>^{-|\alpha|}.
\end{align}
The goal of this section is to prove a Littlewood-Paley type decomposition under a suitable restriction on the initial data. At first we record (without proof) some known results on the pseudodifferential calculus which will be used throughout the paper. We refer to \cite{Robert2,Martinez} for the details of the proof. 
\subsection{Pseudodifferential calculus}
For the metric 
$
g=dx^2/\<x\>^2+d\xi^2/\<\xi\>^2
$ 
and a weight function $m(x,\xi)$ on the phase space $\R^{2d}$, we use H\"ormander's symbol class notation $S(m,g)$, \emph{i.e.}, $a \in S(m,g)$ if and only if $a\in C^\infty(\R^{2d})$ and 
$$
|\dderiv{x}{\xi}{\alpha}{\beta}a(x,\xi)|\le C_{\alpha\beta}m(x,\xi)\<x\>^{-|\alpha|}\<\xi\>^{-|\beta|},\quad \alpha,\beta\in \Z^d_+.
$$
To a symbol $a \in C^\infty(\R^{2d})$ and $h\in(0,1]$, we associate the semiclassical  pseudodifferential operator ($h$-PDO for short) $\Op_h(a)$ defined by
$$
\Op_h(a)f(x)=\frac{1}{(2\pi h)^d}\int e^{i(x-y)\cdot\xi/h}a(x,\xi)f(y)dyd\xi,\quad f\in \S(\R^d).
$$
When $h=1$ we write $\Op(a)=\Op_h(a)$ for simplicity. The Calder\'on-Vaillancourt theorem shows that for any symbol $a\in C^\infty(\R^{2d})$ satisfying
$
|\dderiv{x}{\xi}{\alpha}{\beta}a(x,\xi)|\le C_{\alpha\beta},
$
$\Op_h(a)$ is extended to a bounded operator on $L^2(\R^d)$ uniformly with respect to $h\in(0,1]$.
Moreover, for any symbol $a$ satisfying 
$$
|\dderiv{x}{\xi}{\alpha}{\beta}a(x,\xi)|\le C_{\alpha\beta}\<\xi\>^{-\gamma},\quad \gamma >d,
$$
$\Op_h(a)$ is extended to a bounded operator from $L^q(\R^d)$ to $L^r(\R^d)$ with the following  bounds:
\begin{align}
\label{pseudodifferential_calculus}
\norm{\Op_h(a)}_{L^q \to L^r} \le C_{qr}h^{-d(1/q-1/r)},\quad 1\le q\le r \le \infty,
\end{align}
where $C_{qr}>0$ is independent of $h\in(0,1]$. These bounds follow from the Schur lemma and an interpolation (see, \emph{e.g.},  \cite[Proposition 2.4]{Bouclet_Tzvetkov_1}). 

For two symbols $a\in S(m_1,g)$ and $b\in S(m_2,g)$, the composition $\Op_h(a)\Op_h(b)$ is also a $h$-PDO and written in the form $\Op_h(c)=\Op_h(a)\Op_h(b)$ with a symbol $c\in S(m_1m_2,g)$ given by $c(x,\xi)=e^{ihD_\eta D_z}a(x,\eta)b(z,\xi)|_{z=x,\eta=\xi}$. Moreover, $c(x,\xi)$ has the following expansion
\begin{align}
\label{asymptotic_expansion_composition}
c=\sum_{|\alpha|=0}^{N-1}\frac{h^{|\alpha|}}{i^{|\alpha|}\alpha!}\partial_\xi^\alpha a\cdot \partial_x^\alpha b+h^Nr_N\ \text{with}\ r_N\in S(\<x\>^{-N}\<\xi\>^{-N}m_1m_2,g).
\end{align}
The symbol of the adjoint $\Op_h(a)^*$ is given by $a^*(x,\xi)=e^{ihD_\eta D_z}a(z,\eta)|_{z=x,\eta=\xi}\in S(m_1,g)$ which has the expansion
\begin{align}
\label{asymptotic_expansion_adjoint}
a^*=\sum_{|\alpha|=0}^{N-1}\frac{h^{|\alpha|}}{i^{|\alpha|}\alpha!}\partial_\xi^\alpha\partial_x^\alpha a+h^Nr_N^*\ \text{with}\ r_N^*\in S(\<x\>^{-N}\<\xi\>^{-N}m_1,g).
\end{align}

\subsection{Littlewood-Paley decomposition}
As we mentioned in the outline of the paper, $H$ is not bounded below in general and we hence cannot expect that the Littlewood-Paley decomposition associated with $H$, which is of the form
$$
\norm{u}_{L^q} \le C_q\norm{u}_{L^2}+C_q\bigg(\sum_{j=0}^\infty \norm{f(2^{-2j}H)u}_{L^q}^2\bigg)^{1/2},
$$
hold if $q\neq2$. The standard Littlewood-Paley decomposition associated with $H_0$ also does not work well in our case, since the commutator of $H$ with the Littlewood-Paley projection $f(2^{-2j}H_0)$ can grow at spatial infinity. To overcome this difficulty, let us introduce an additional localization as follows. Given a parameter $\ep>0$ and a cut-off function $\varphi \in C^\infty_0(\R_+)$ such that $\varphi \equiv 1$ on $[0,1/2]$ and $\supp \varphi\subset[0,1]$, we define $\psi_\ep(x,\xi)$ by
$$
\psi_\ep(x,\xi)=\varphi\left(\frac{\<x\>}{\ep|\xi|}\right).
$$
 It is easy to see that, for each $\ep>0$, $\psi_\ep\in S(1,g)$ and is supported in $\{(x,\xi)\in\R^{2d};\<x\><\ep|\xi|\}$. Moreover, for sufficiently small $\ep>0$, $p(x,\xi)$ is uniformly elliptic on the support of $\psi_\ep$ and $\Op(\psi_\ep)H$ thus is essentially bounded below. 
 
 In this subsection we prove a Littlewood-Paley type decomposition on the range of $\Op(\psi_\ep)$. We begin with the following proposition which tells us that, for any $f\in C_0^\infty(\R)$ and $h\in(0,1]$, $\Op(\psi_\ep)f(h^2H)$ is well approximated in terms of the $h$-PDO.

%proposition
\begin{proposition}					
\label{proposition_functional_calculus}
There exists $\ep>0$ such that, for any $f\in C_0^\infty(\R)$ with $\supp f \Subset (0,\infty)$, we can  construct bounded families $\{a_{h,j}\}_{h\in(0,1]} \subset \bigcap_{M\ge0}S(\<x\>^{-j}\<\xi\>^{-M},g)$, $j\ge0$, such that\\
\emph{(1)} $a_{h,0}$ is given explicitly by $a_{h,0}(x,\xi)=\psi_{\ep}(x,\xi/h)f(p_h(x,\xi))$. Moreover, 
\begin{align*}
\supp a_{h,j}
\subset \supp\psi_{\ep}(\cdot,\cdot/h) \cap\supp f(p_h)
\subset\{(x,\xi)\in\R^{2d};\<x\><1/h,\ |\xi|\in I\},
\end{align*}
for some relatively compact open interval $I\Subset(0,\infty)$. In particular, we have
$$
\norm{\Op_h(a_{h,j})}_{L^{q'}\to L^q} \le C_{jqq'}h^{-d(1/q'-1/q)},\ 1\le q'\le q\le \infty,
$$
uniformly in $h\in (0,1]$. \\
\emph{(2)} For any integer $N>d+2$, we set $a_h=\sum_{j=0}^{N-1}h^ja_{h,j}$. Then,
$$
\norm{\Op(\psi_\ep) f(h^2H)-\Op_h(a_h)}_{L^2\to L^q} \le C_{qN}h^{2},\quad 2\le q\le \infty,
$$
uniformly in $h\in (0,1]$. 
\end{proposition} 

The following is an immediate consequence of this proposition.
\begin{corollary}
\label{corollary_functional_calculus}
For any $2\le q\le\infty$ and $h\in(0,1]$, $\Op(\psi_\ep)f(h^2H)$ is bounded from $L^2(\R^d)$ to $L^q(\R^d)$ and satisfies
$$
\norm{\Op(\psi_\ep)f(h^2H)}_{L^2\to L^q} \le C_{q}h^{-d(1/2-1/q)},
$$
where $C_q>0$ is independent of $h\in(0,1]$. 
\end{corollary}

For the low energy part we have the following
\begin{lemma}							%{lemma}
\label{lemma_functional_calculus}
For any $f_0\in C_0^\infty(\R)$ and $2\le q\le\infty$, we have
$$
\norm{\Op(\psi_\ep)f_0(H)}_{L^2\to L^q} \le C_{q}.
$$
\end{lemma}

\begin{remark}
If $V,A\equiv0$, then Proposition \ref{proposition_functional_calculus}, Corollary \ref{corollary_functional_calculus} and Lemma \ref{lemma_functional_calculus} hold without the additional term $\Op(\psi_\ep)$. Moreover, in this case we see that the remainder satisfies 
$$
\norm{f(h^2H)-\Op_h(a_h)}_{L^2\to L^q} \le C_{qN}h^{N-d(1/2-1/q)}.
$$
We refer to  \cite{BGT} (for the case on compact manifolds without boundary) and to \cite{Bouclet_Tzvetkov_1} (for the case with metric perturbations on $\R^d$). For more general cases with Laplace-Beltrami operators on non-compact manifolds with ends, we refer to \cite{Bouclet3,Bouclet2}. Because of this result, we believe Proposition \ref{proposition_functional_calculus} is far from sharp. However, the bounds
$$
\norm{\Op(\psi_\ep) f(h^2H)-\Op_h(a_h)}_{L^2\to L^q} \le C_{qN}h,\ 2\le q\le\infty,
$$
are sufficient to obtain our Littlewood-Paley type decomposition (Proposition \ref{proposition_Littlewood_Paley}). For more details, we refer to Burq-G\'erard-Tvzetkov \cite[Corollary 2.3]{BGT}. 
\end{remark}

\begin{proof}[Proof of Proposition \ref{proposition_functional_calculus}]
We write 
$$
\Op(\psi_\ep)=\Op_h(\psi_{\ep/h}),\quad h\in(0,1],
$$
where $\psi_{\ep/h}(x,\xi)=\psi_\ep(x,\xi/h)$ satisfies $\supp \psi_{\ep/h}\subset\{h\<x\><\ep|\xi|\}$ and 
\begin{align}
\label{proof_proposition_functional_calculus_1}
|\dderiv{x}{\xi}{\alpha}{\beta}\psi_{\ep/h}(x,\xi)|
\le C_{\alpha\beta\ep} h^{-|\beta|}\<x\>^{-|\alpha|}\<\xi/h\>^{-|\beta|}
\le C_{\alpha\beta\ep} \<x\>^{-|\alpha|}(h+|\xi|)^{-|\beta|}.
\end{align}
We have by using the Helffer-Sj\"ostrand formula (see, Helffer-Sj\"ostrand \cite{Helffer_Sjostrand}) that
$$
\Op_h(\psi_{\ep/h})f(h^2H)=-\frac{1}{2\pi i}\int_{\C}\frac{\partial{\wtilde{f}}}{\partial \bar{z}}(z)\Op_h(\psi_{\ep/h})(h^2H-z)^{-1}dz\wedge d\bar{z},
$$
where $\wtilde{f}(z)$ is an almost analytic extension of $f(\lambda)$. Since $ f\in C^\infty_0(\R)$, $\wtilde{f}(z)$ is also compactly supported and satisfies 
$$
\partial_{\bar{z}}\wtilde{f}(z) =O(|\Im z|^M)
$$
 for any $M>0$. We may assume $|z|\le C$ on $\supp \wtilde{f}$ with some $C>0$. In order to use this formula we shall construct a semiclassical approximation of $\Op_h(\psi_{\ep/h})(h^2H-z)^{-1}$, in terms of the $h$-PDO, for $z\in \C\setminus[0,\infty)$ with $|z|\le C$. Although the method is based on the standard semiclassical parametrix construction (see, \emph{e.g.}, \cite{Robert2,BGT}), we give the detail of the proof since $\psi_{\ep/h}$ is not uniformly bounded in $S(1,g)$ with respect to $h\in(0,1]$. 

We first study the symbol of the resolvent $(h^2H-z)^{-1}$. Let $p_h$ and $p_{1,h}$ be as in \eqref{symbols_2} so that $h^2H=\Op_h(p_h)+h\Op_h(p_{1,h})$. Since 
$$
h|A(x)|\lesssim |\xi|,\ h^2|V(x)|\lesssim |\xi|^2,
$$
on $\supp \psi_{\ep/h}$, we obtain by \eqref{estimates_on_symbols_0} that
\begin{align}
\label{proof_proposition_functional_calculus_2}
|\dderiv{x}{\xi}{\alpha}{\beta}p_{h}(x,\xi)|
&\le C_{\alpha\beta}\<x\>^{-|\alpha|}|\xi|^{2-|\beta|}\quad\text{if}\ |\beta|\le2,\\
\label{proof_proposition_functional_calculus_3}
|\dderiv{x}{\xi}{\alpha}{\beta}p_{1,h}(x,\xi)|
&\le C_{\alpha\beta}\<x\>^{-1-|\alpha|}|\xi|^{1-|\beta|}\quad\text{if}\ |\beta|\le1,
\end{align}
uniformly in $(x,\xi)\in \supp\psi_{\ep/h}$ and $h\in(0,1]$. Moreover, if $\ep>0$ is sufficiently small then the uniform ellipticity of $k$ implies that $p_h$ is also uniformly elliptic on $\supp\psi_{\ep/h}$: 
$$
C_1^{-2}|\xi|^2 \le p_h(x,\xi)\le C_1^2 |\xi|^2\quad\text{if}\ h\<x\><\ep|\xi|,
$$
with some $C_1>0$, which particularly implies 
\begin{align}
\label{resolvent}
\frac{1}{|p_h(x,\xi)-z|}\lesssim
\begin{cases}
 |\Im z|^{-1}&\text{if}\ |\xi| \le 2C_2,\\
 \<\xi\>^{-2}&\text{if}\ |\xi| \ge 2C_2
\end{cases}
\end{align}
for $(x,\xi)\in\supp \psi_{\ep/h}$, $z\notin \R$ and $|z|\le C$, with some $C_2>0$. 

Let us now consider a sequence of symbols $q_{j}^h=q_{j}^h(z,x,\xi)$ (depending holomorphically on $z\notin\R$) defined inductively by
\begin{align*}
q_{0}^h&=\frac{\psi_{\ep/h}}{p_h-z},\\
q_{1}^h&=-\frac{1}{p_h-z}\bigg(\sum_{|\alpha|=1}i^{-1}\partial_\xi^\alpha q_{0}^h\cdot\partial_x^\alpha p_h+q_{0}^h\cdot p_{1,h}\bigg),\\
q_{j}^h&=-\frac{1}{p_h-z}\bigg(\sum_{|\alpha|+k=j,\atop|\alpha|\ge1}\frac{i^{-|\alpha|}}{\alpha!}\partial_\xi^\alpha q_{k}^h\cdot\partial^\alpha_x p_h+\sum_{|\alpha|+k={j-1},}\frac{i^{-|\alpha|}}{\alpha!}\partial_\xi^\alpha q_{k}^h\cdot\partial^\alpha_x p_{1,h}\bigg),\ j\ge2.
\end{align*}
We then learn by \eqref{proof_proposition_functional_calculus_1}, \eqref{proof_proposition_functional_calculus_2} and \eqref{resolvent} that
\begin{align}
\label{q_0_1}
|\dderiv{x}{\xi}{\alpha}{\beta}q_{0}^h(z,x,\xi)|
&\le C_{\alpha\beta\ep}
\begin{cases}
\<x\>^{-|\alpha|}(h+|\xi|)^{-|\beta|}|\Im z|^{-1-|\alpha+\beta|}&\text{if}\ |\xi| \le 2C_2,\\
\<x\>^{-|\alpha|}\<\xi\>^{-|\beta|-2}&\text{if}\ |\xi| \ge 2C_2,
\end{cases}\\
\nonumber
&\le C_{\alpha\beta \ep} 
\<x\>^{-|\alpha|}(h+|\xi|)^{-|\beta|}|\Im z|^{-1-|\alpha+\beta|}
\end{align}
for $z\notin\R$ with $|z|\le C$ and $h\in(0,1]$. We similarly obtain by using \eqref{proof_proposition_functional_calculus_2}, \eqref{proof_proposition_functional_calculus_3} and \eqref{q_0_1} that if $h|\xi|\le 2C_2$ then
\begin{align*}
|\dderiv{x}{\xi}{\alpha}{\beta}q_{1}^h(z,x,\xi)|
&\le C_{\alpha\beta\ep}
\bigg(\<x\>^{-1-|\alpha|}(h+|\xi|)^{-1-|\beta|}|\xi|^2|\Im z|^{-3-|\alpha+\beta|}\\
&\ \ \ \ \ +
\<x\>^{-1-|\alpha|}(h+|\xi|)^{-|\beta|}(h+|\xi|)|\Im z|^{-2-|\alpha+\beta|}\bigg)\\
&\le C_{\alpha\beta\ep}
(h+|\xi|)^2\<x\>^{-1-|\alpha|}(h+|\xi|)^{-1-|\beta|}|\Im z|^{-3-|\alpha+\beta|},
\end{align*}
for $z\notin\R$ with $|z|\le C$ and $h\in(0,1]$. Here note that, in this case, $(h+|\xi|)^{-1}$ may have a singularity at $\xi=0$ as $h\to+0$. In order to prove the remainder estimate, we will remove this singularity by using a rescaling $\xi\mapsto h\xi$ (see the estimates \eqref{proposition_functional_calculus_2}). For $h|\xi|\ge 2C_2$, $q_{1}^h$ does not have such a singularity and satisfies
\begin{align*}
|\dderiv{x}{\xi}{\alpha}{\beta}q_{1}^h(z,x,\xi)|
\le C_{\alpha\beta\ep}
\<x\>^{-1-|\alpha|}\<\xi\>^{-|\beta|-4}|\xi|
\le C_{\alpha\beta\ep}
\<x\>^{-1-|\alpha|}\<\xi\>^{-|\beta|-3}.
\end{align*}
uniformly in $z\notin\R$ with $|z|\le C$ and $h\in(0,1]$. 
Since $1\lesssim h+|\xi|$ if $h|\xi| \gtrsim1$, summarizing these we have
$$
|\dderiv{x}{\xi}{\alpha}{\beta}q_{1}^h(z,x,\xi)|
\le C_{\alpha\beta\ep}
\<x\>^{-1-|\alpha|}(h+|\xi|)^{1-|\beta|}|\Im z|^{-3-|\alpha+\beta|},\quad 
z\notin\R,\ |z|\le C,\ h\in(0,1].
$$
The estimates \eqref{q_0_1} and a direct computation also show that $q_{1}^h$ is of the form
$$
q_{1}^h=q_{11}^h(p_h-z)^{-3}+q_{10}^h(p_h-z)^{-2},
$$
where $q_{1k}^h$ are supported in $\supp \psi_{\ep/h}$, independent of $z$ and satisfy 
$$
|\dderiv{x}{\xi}{\alpha}{\beta}q_{1k}^h(x,\xi)|
\le C_{\alpha\beta \ep} 
\<x\>^{-1-|\alpha|}(h+|\xi|)^{-|\beta|}\<\xi\>^{N_1(k)},\quad h\in(0,1],
$$
with some positive integer $N_1(k)>0$. For $j\ge2$, an induction argument yields that
\begin{align}
\label{q_j_1}
|\dderiv{x}{\xi}{\alpha}{\beta}q_{j}^h(z,x,\xi)|\le C_{\alpha\beta\ep}
\<x\>^{-j-|\alpha|}(h+|\xi|)^{2-j-|\beta|}|\Im z|^{-2j-1-|\alpha+\beta|},\quad j\ge2,
\end{align}
for $z\notin\R$ with $|z|\le C$ and $h\in(0,1]$. It also follows from an induction on $j$ that there exist a sequence of $z$-independent symbols $(q_{jk}^h)_{k=0}^{j}$ supported in $\supp \psi_{\ep/h}$ and satisfying
\begin{align}
\label{q_j_2}
|\dderiv{x}{\xi}{\alpha}{\beta}q_{jk}^h(x,\xi)|
\le C_{\alpha\beta\ep}
\<x\>^{-j-|\alpha|}(h+|\xi|)^{-|\beta|}\<\xi\>^{N_j(k)}
\end{align}
with some $N_j(k)>0$, such that $q_{j}^h$ is of the form
$$
q_{j}^h=\sum_{k=0}^{j}q_{jk}^h(p_h-z)^{-j-k-1}.
$$

Rescaling $\xi\mapsto h\xi$, we learn by \eqref{q_0_1} and \eqref{q_j_1} that 
$$
q_0^h(z,x,h\xi)\in S(1,g),\quad h^jq_j^h(z,x,h\xi)\in S(h^2\<x\>^{-j}\<\xi\>^{2-j},g),
$$
with uniform bounds in $h$ and polynomially bounds in $|\Im z|^{-1}$. Then, by the construction of $q_j^h$, the standard symbolic calculus (not in the semiclassical regime) and the fact that
$$
\Op(h^jq_j^h(z,x,h\xi))=h^j\Op_h(q_j^h),
$$
we obtain
$$
\Op(\psi_\ep)=\sum_{j=0}^{N-1}h^j\Op_h(q_{j}^h)(h^2H-z)+h^2\Op(r_{h,N,z}),\quad N\ge1,
$$
with some $r_{h,N,z}\in S(\<x\>^{-N}\<\xi\>^{2-N},g)$ satisfying
\begin{align}
\label{proposition_functional_calculus_2}
|\dderiv{x}{\xi}{\alpha}{\beta}r_{h,N,z}(x,\xi)|\le C_{\alpha\beta\ep N}
\<x\>^{-N-|\alpha|}\<\xi\>^{2-N-|\beta|}|\Im z|^{-2N-1-|\alpha+\beta|},
\end{align}
where $C_{\alpha\beta\ep N}>0$ may be taken uniformly in $h\in(0,1]$, $z\in \C\setminus\R$ with $|z|\le C$ and $x,\xi\in\R^d$. 

We now use the Helffer-Sj\"ostrand formula to obtain
$$
\Op(\psi_\ep) f(h^2H)=\sum_{j=0}^{N-1}h^j\Op_h(a_{h,j})+h^2R(h,N),
$$
where 
\begin{align*}
a_{h,0}(x,\xi)&=\psi_{\ep/h}(x,\xi)(f\circ p_h)(x,\xi),\\
a_{h,j}(x,\xi)&=\sum_{k=0}^{j}\frac{(-1)^{k+j}}{(k+j)!}q_{jk}^h(x,\xi)(f^{(j+k)}\circ p_h)(x,\xi),\ 1\le j\le N-1,\\
R(h,N)&=-\frac{1}{2\pi i}\int_{\C}\frac{\partial{\wtilde{f}}}{\partial \bar{z}}(z)\Op_h(r_{h,N,z})(h^2H-z)^{-1}dz\wedge d\bar{z}. 
\end{align*}
Since $\supp q_{jk}\subset \supp \psi_{\ep/h}\subset \{h\<x\><\ep|\xi|\}$ and $p_h$ is uniformly elliptic (\emph{i.e.}, $p_h\approx |\xi|^2$) on the latter region, taking $\ep>0$ smaller if necessary we have
$$
a_{h,j}\subset \supp \psi_{\ep/h} \cap\supp f(p_h)\subset\{(x,\xi);|x|<1/h,\ C_0^{-1}\le |\xi|\le C_0\}
$$ 
with some positive constant $C_0>0$, which, combined with \eqref{q_j_2}, particularly implies $q_{jk}^h\in S(\<x\>^j,g)$ since $h+|\xi|\gtrsim \<\xi\>$ on $\supp \psi_{\ep/h} \cap\supp f(p_h)$. $\{a_{h,j}\}_{h\in(0,1]}$ is hence bounded in $\bigcap_{M\ge0}S(\<x\>^{-j}\<\xi\>^{-M},g)$. By virtue of \eqref{pseudodifferential_calculus}, we also obtain
$$
\norm{\Op_h(a_{h,j})}_{L^{q'} \to L^q} \le C_{jqq'} h^{-d(1/q'-1/q)},\quad h\in(0,1],\ 1\le q'\le q\le \infty.
$$

Finally, we shall prove the estimate on the remainder $R(h,N)$. If we choose $N>d+2$, then \eqref{proposition_functional_calculus_2} and \eqref{pseudodifferential_calculus} (with $h=1$) imply
$$
\norm{\Op(r_{h,N,z})}_{L^2 \to L^q} \le C_{qN}|\Im z|^{-n(N,q)}, \quad 2\le q \le \infty,
$$
with some positive integer $n(N,q)\ge 2N+1$, where $C_{qN}>0$ is independent of $h$. Using the bounds $\norm{(h^2H-z)^{-1}}_{L^2 \to L^2} \le |\Im z|^{-1}$, $|\partial_{\bar{z}}\wtilde{f}(z)|\le C_M|\Im z|^M$ for any $M\ge0$ and the fact that $\wtilde{f}$ is compactly supported, we conclude that
\begin{align*}
\norm{R(h,N)}_{L^2 \to L^q}
&\le C_M\int_{\supp\wtilde{f}} |\Im z|^M\norm{\Op(r_{h,N,z})}_{L^2 \to L^q}\norm{(h^2H-z)^{-1}}_{L^2 \to L^2}dz\wedge d\bar{z}\\
&\le C_{MNq}\int_{\supp\wtilde{f}} |\Im z|^{M-n(N,q)-1}dz\wedge d\bar{z}\\
&\le C_{MNq},
\end{align*}
provided that $M$ is large enough. We complete the proof. 
\end{proof}

\begin{proof}[Proof of Lemma \ref{lemma_functional_calculus}]
By the same argument as above with $h=1$, we can see that
$$
\Op(\psi_\ep)f_0(H)=\sum_{j=0}^{N-1}\Op(a_j)+R(N)
$$
where $a_j\in \bigcap_{M\ge0} S(\<x\>^{-j}\<\xi\>^{-M},g)$ are supported in 
$$
\supp \psi_\ep\cap\supp f_0(p)\subset \{(x,\xi)\in\R^{2d}; \<x\><\ep|\xi|,\ |\xi|\lesssim 1\}
$$
and $R(N)$ satisfies
$$
\norm{R(h,N)}_{L^2 \to L^q}\le
C_{Nq},\quad 2\le q \le\infty,
$$ if $N>d+2$. The assertion then follows from \eqref{pseudodifferential_calculus}. 
\end{proof}

Consider a $4$-adic partition of unity:
$$
f_0(\lambda)+\sum_h f(h^2\lambda)=1,\quad \lambda\in\R,
$$
where $f_{0},f\in C^\infty_0(\R)$ with $\supp f_{0}\subset [-1,1]$, $\supp f\subset[1/4,4]$ and $\sum_h$ means that, in the sum, $h$ takes all negative powers of 2 as values, \emph{i.e.}, 
$\sum_h=\sum_{h=2^{-j},j\ge0}$. Let $F\in C_0^\infty(\R)$ be such that $\supp F \subset [1/8,8]$ and $F\equiv1\ \text{on}\ \supp f$. 
The spectral decomposition theorem implies
$$
1=f_0(H)+\sum_h f(h^2H)=f_0(H)+\sum_h F(h^2H)f(h^2H).
$$
Let $a_h \in S(1,g)$ be as in Proposition \ref{proposition_functional_calculus} with $f=F$. Using Proposition \ref{proposition_functional_calculus}, we obtain a Littlewood-Paley type estimates on a range of $\Op(\psi_\ep)$. 

%proposition
\begin{proposition}					
\label{proposition_Littlewood_Paley}
For any $2\le q<\infty$, 
$$
\norm{\Op(\psi_\ep)u}_{L^q(\R^d)} \le C_q\norm{u}_{L^2(\R^d)}+C_q\bigg(\sum_{h}\norm{\Op_h(a_h)f(h^2 H)u}_{L^q(\R^d)}\bigg)^{1/2}.
$$
\end{proposition}

\begin{proof}
The proof is same as that of \cite[Corollary 2.3]{BGT} and we omit details.
\end{proof}

%corollary
\begin{corollary}
\label{corollary_Littlewood_Paley}
Let $\ep>0$ and $\psi_\ep$ be as above and $\chi_\ep=1-\psi_\ep$. Let $\rho\in C^\infty(\R^d)$ be such that 
$$
|\partial_x^\alpha \rho(x)|\le C_\alpha\<x\>^{-|\alpha|},\quad \alpha\in \Z^d_+.
$$
Then, for any $T>0$ and any $(p,q)$ satisfying $p\ge2,q<\infty$ and $2/p=d(1/2-1/q)$, there exists $C_T>0$ such that
\begin{align*}
\norm{\rho e^{-itH}\varphi}_{L^p([-T,T];L^q(\R^d))}
&\le C_T\norm{\varphi}_{L^2(\R^d)}+C\norm{\Op(\chi_\ep)e^{-itH}\varphi}_{L^p([-T,T];L^q(\R^d))}\\
&+C\bigg(\sum_{h}\norm{\Op_h(a_h)e^{-itH} f(h^2 H)\varphi}_{L^p([-T,T];L^q(\R^d))}^2\bigg)^{1/2}, 
\end{align*}
where $a_h$ is given by Proposition \ref{proposition_functional_calculus} with $\psi_\ep$ replaced by $\rho \psi_\ep$. In particular, $a_h(x,\xi)$ is supported in $\supp \rho(x) \psi(x,\xi/h)F(p_h(x,\xi))$. 
\end{corollary}

\begin{proof}
This proposition follows from the $L^2$-boundedness of $e^{-itH}$, Propositions \ref{proposition_functional_calculus} and \ref{proposition_Littlewood_Paley} (with $\psi_\ep$ replaced by $\rho \psi_\ep$) and the Minkowski inequality.
\end{proof}

%%%%%%%%%%%%%%%%%%%%%%%%%%%%%%%%%%%%%%%%%%%%%%%%      Proof of Theorem \ref{theorem_1}                              %%%%%%%%%                                             %%%%%%%%%%%%%%%%%%%%%%%%%%%%%%%%%%%%%%%%%%%

\section{Proof of Theorem \ref{theorem_1}}
\label{section_proof_theorem_1}
In this section we prove Theorem \ref{theorem_1} under Assumption \ref{assumption_A} with $\mu>0$. We first state two key estimates which we will prove in later sections. For $R>0$, an open interval $I\Subset(0,\infty)$ and $\sigma\in(-1,1)$, we define the outgoing and incoming regions $\Gamma^\pm(R,I,\sigma)$ by 
$$
\Gamma^\pm(R,I,\sigma):=\left\{(x,\xi)\in\R^{2d};|x|>R,\ |\xi|\in I,\  \pm\frac{x\cdot\xi}{|x||\xi|}>-\sigma\right\}, 
$$ 
respectively. We then have the following (local-in-time) dispersive estimates:

%theorem
\begin{proposition}
\label{theorem_dispersive_IK}
Suppose that $H$ satisfies Assumption \ref{assumption_A} with $\mu>0$. Let $I\Subset(0,\infty)$ and $\sigma\in(-1,1)$. Then, for sufficiently large $R\ge1$, small $h_0>0$ and any symbols $a^\pm_h\in S(1,g)$ supported in $\Gamma^\pm(R,I,\sigma)\cap \{x;|x|<1/h\}$, we have
$$
\norm{\Op_h(a_h^\pm)e^{-itH}\Op_h(a_h^\pm)^*}_{L^1\to L^\infty} \le C|t|^{-d/2},\quad 0<  |t| \le 1,
$$
uniformly with respect to $h\in(0,h_0]$.
\end{proposition}

We prove this proposition in Section \ref{section_IK}. In the region $\{|x|\gtrsim|\xi|\}$, we have the following (short-time) dispersive estimates:
\begin{proposition}
\label{theorem_dispersive_WKB}
Suppose that $H$ satisfies Assumption \ref{assumption_A} with $\mu\ge0$. Let us fix arbitrarily $\ep>0$. Then, there exists $t_\ep>0$ such that, for any symbol $\chi_\ep\in S(1,g)$ supported in 
$
\{(x,\xi); \<x\>\ge\ep|\xi|\}
$, 
we have
$$
\norm{\Op(\chi_\ep)e^{-itH}\Op(\chi_\ep)^*}_{L^1\to L^\infty} \le C_\ep|t|^{-d/2},\quad 0<  |t| \le t_\ep.
$$
\end{proposition}

We prove this proposition in Section \ref{section_WKB}. 

\begin{proof}[Proof of Theorem \ref{theorem_1}]
Taking $\rho\in C^\infty(\R^d)$ so that $0\le \rho(x)\le1$, $\rho(x)=1$ for $|x|\ge1$ and $\rho(x)=0$ for $|x|\le1/2$, we set $\rho_R(x)=\rho(x/R)$. 
In order to prove Theorem \ref{theorem_1}, it suffices to show
$$
\norm{\rho_Re^{-itH}\varphi}_{L^p([-T,T];L^q(\R^d))} \le C_T\norm{\varphi}_{L^2(\R^d)},
$$
for sufficiently large $R\ge1$. We also may assume without loss of generality that $T>0$ is sufficiently small. Indeed, if the above estimate holds on $[-T_0,T_0]$ with some $T_0>0$ then we obtain by the unitarity of $e^{-itH}$ on $L^2$ that, for any $T>T_0$,
\begin{align*}
\norm{\rho_Re^{-itH}\varphi}_{L^p([-T,T];L^q(\R^d))}^p
&\lesssim \sum_{k=-[T/T_0]}^{[T/T_0]+1}\norm{\rho_Re^{-itH}e^{-i(k+1)H}\varphi}_{L^p([-T_0,T_0];L^q(\R^d))}^p\\
&\lesssim (T/T_0)C_{T_0}^p\norm{\varphi}_{L^2(\R^d)}^p.
\end{align*}
Let $a_h$ be as in Proposition \ref{proposition_functional_calculus}. Replacing $\psi_\ep$ with $\rho_R\psi_\ep$ and taking $\ep>0$ smaller if necessary, we may assume without loss of generality that 
$
\supp a_h \subset  \{(x,\xi);R<|x|<1/h, |\xi|\in I\}
$ for some open interval $I\Subset(0,\infty)$. Choosing $\theta^\pm \in C^\infty([-1,1])$ so that $\theta^++\theta^-=1$, $\theta^+=1\ \text{on}\ [1/2,1]$ and $\theta^+=0\ \text{on}\ [-1,-1/2]$, we set 
$a_h^\pm(x,\xi)=a_h(x,\xi)\theta^\pm(\hat{x}\cdot\hat{\xi})$, where $\hat{x}=x/|x|$. It is clear that $\{a_h^\pm\}_{h\in(0,1]}$ is bounded in $S(1,g)$ and $\supp a_h^\pm \subset \Gamma^\pm(R,I,1/2)\cap \{x;|x|<1/h\}$, and that $a_h=a_h^++a_h^-$. We now apply Proposition \ref{theorem_dispersive_IK} to $a^\pm_h$ and obtain the local-in-time dispersive estimate for $\Op_h(a_h^\pm)e^{-itH}\Op_h(a_h^\pm)^*$ (uniformly in $h\in (0,h_0]$), which, combined with the $L^2$-boundedness of $\Op_h(a^\pm_h)e^{-itH}$ and the abstract Theorem due to Keel-Tao \cite{Keel_Tao}, implies Strichartz estimates for $\Op_h(a_h)e^{-itH}$:
\begin{align*}
\norm{\Op_h(a_h)e^{-itH}\varphi}_{L^p([-1,1];L^q(\R^d))}
&\le \sum_\pm \norm{\Op_h(a_h^\pm)e^{-itH}\varphi}_{L^p([-1,1];L^q(\R^d))}\\
&\le C\norm{\varphi}_{L^2(\R^d)},
\end{align*}
uniformly with respect to $h\in(0,h_0]$. Since $\Op_h(a_h)$ is bounded from $L^2(\R^d)$ to $L^q(\R^d)$ with the bound of order $O(h^{-d(1/2-1/q)})$, for $h_0<h\le1$ we have
$$
\sum_{h_0< h\le1}\norm{\Op_h(a_h)e^{-itH} f(h^2 H)\varphi}_{L^p([-1,1];L^q(\R^d))}^2 \le C(h_0)\norm{\varphi}_{L^2(\R^d)}^2. 
$$
with some $C(h_0)>0$. Using these two bounds, we obtain
\begin{align*}
&\sum_{h}\norm{\Op_h(a_h)e^{-itH} f(h^2 H)\varphi}_{L^p([-1,1];L^q(\R^d))}^2\\
&\le C\sum_{0<h<h_0}\norm{f(h^2 H)\varphi}_{L^2(\R^d))}^2+C(h_0)\norm{\varphi}_{L^2(\R^d)}^2\\
&\le C\norm{\varphi}_{L^2(\R^d)}^2.
\end{align*}
On the other hand, Strichartz estimates for $\Op(\chi_\ep)e^{-itH}$ is an immediate consequence of Proposition \ref{theorem_dispersive_WKB}. By virtue of Corollary \ref{corollary_Littlewood_Paley}, we complete the proof.
\end{proof}

%%%%%%%%%%%%%%%%%%%%%%%%%%%%%%%%%%%%%%%%%%%%%%%%      Parametrix of Isozaki-Kitada type	                                                                             %%%%%%%%%%%%%%%%%%%%%%%%%%%%%%%%%%%%%%%%%%%%%%

\section{Semiclassical approximations for outgoing propagators}
\label{section_IK}
Throughout this section we assume Assumption \ref{assumption_A} with $\mu>0$. We here study the behavior of $e^{-itH}\Op_h(a_h^\pm)^*$, where $a_h^\pm\in S(1,g)$ are supported in $\Gamma^\pm(R,I,\sigma)\cap\{|x|<1/h\}$, respectively. The main goal of this section is to prove Proposition \ref{theorem_dispersive_IK}. For simplicity, we consider the outgoing propagator $e^{-itH}\Op_h(a_h^+)^*$ for $0\le t \le 1$ only, and the proof for the incoming case is analogous. 

In order to prove dispersive estimates, we construct a semiclassical approximation for the outgoing propagator $e^{-itH}\Op_h(a_h^+)^*$ by using the method of Isozaki-Kitada. Namely, rescaling $t\mapsto th$ and setting $H^h=h^2H$, $H_0^h=-h^2\Delta/2$, we consider an approximation for the semiclassical propagator $e^{-itH^h/h}\Op_h(a_h^+)^*$ of the following form
$$
e^{-itH^h/h}\Op_h(a_h^+)^*=J_h(S^+_h,b^+_h)e^{-itH_0^h/h}J_h(S^+_h,c^+_h)^*+O(h^N),\quad 0\le t\le h^{-1},
$$
where $S^+_h$ solves a suitable Eikonal equation in the outgoing region and $J(S^+_h,w)$ is the corresponding semiclassical Fourier integral operator ($h$-FIO for short):
$$
J_h(S^+_h,w)f(x)=(2\pi h)^{-d}\int e^{i(S^+_h(x,\xi)-y\cdot\xi)/h}w(x,\xi)f(y)dyd\xi. 
$$
Such approximations (uniformly in time) have been studied by \cite{Robert_Tamura} for  Schr\"odinger operators with long-range potentials, and by \cite{Robert2,Robert,Bouclet_Tzvetkov_1} for the case of long-range metric perturbations. We also refer to the original paper by Isozaki-Kitada \cite{Isozaki_Kitada} in which the existence and asymptotic completeness of modified wave operators (with time-independent modifiers) were established for the case of Schr\"odinger operators with long-range potentials. We note that, in these cases, we do not need the additional restriction of the initial data in $\{|x|<1/h\}$. 
On the other hand, in a recent paper \cite{Mizutani2}, we constructed such approximations (locally in time) for the case with long-range metric perturbations, combined with potentials growing subquadratically at infinity, under the additional restriction on the initial data into $\{|x|<1/h\}$. 

As we mentioned in the outline of the paper, we first construct an approximation for the modified propagator $e^{-it\wtilde{H}^h/h}$, where $\wtilde{H}^h$ is defined as follows. Taking arbitrarily a cut-off function $\psi \in C^\infty_0(\R^d)$ such that $0\le\psi\le1$, $\psi\equiv 1\ \text{for}\ |x|\le1/2$ and $\psi\equiv0\ \text{for}\ |x|\ge1$,  
we define truncated electric and magnetic potentials, $V_{h}$ and $A_{h}=(A_{h,j})_j$ by
$V_{h}(x):=\psi(hx/L)V(x)$, $A_{h,j}(x)=\psi(hx/L)A_j(x)$, 
respectively. It is easy to see that 
$$
V_{h}\equiv V,\ A_{h,j}\equiv A_j\ \text{on}\ \{|x|\le L/(2h)\},\quad 
\supp A_{h,j},\ \supp V_{h}\subset\{|x|\le L/h\},
$$
and that, for any $\alpha\in\Z^d_+$ there exists $C_{L,\alpha}>0$, independent of $x,h$, such that
\begin{align}
\label{estimates_on_V_h}
h^2|\partial_x^\alpha V_{h}(x)|+h|\partial_x^\alpha A_{h}(x)|&\le C_{\alpha,L} \<x\>^{-\mu-|\alpha|}.
\end{align}
Let us define $\wtilde{H}^h$ by
$$
\wtilde{H}^h=\frac 12\sum_{j,k=1}^d (-ih\partial_j-hA_{h,j}(x)) g^{jk}(x) (-ih\partial_k-hA_{h,k}(x))+h^2V_{h}(x).
$$
We consider $\wtilde{H}^h$ as a ``semiclassical" Schr\"odinger operator with $h$-dependent electromagnetic potentials $h^2V_{h}$ and $hA_{h}$. By virtue of the estimates on $g^{jk},A_{h}$ and $V_{h}$, $\wtilde{H}^h$ can be regarded as a long-range perturbation of the semiclassical free Schr\"odinger operator $H_0^h=-h^2\Delta/2$. Such a type modification has been used to prove Strichartz estimates and local smoothing effects (with loss of derivatives) for Schr\"odinger equations with superquadratic potentials (see, Yajima-Zhang \cite[Section 4]{Yajima_Zhang}). Let us denote by $\wtilde{p}_{h}$ the corresponding modified symbol:
\begin{align}
\label{symbol1}
\wtilde{p}_{h}(x,\xi)=\frac12\sum_{j,k=1}^dg^{jk}(x)(\xi_j-hA_{h,j}(x))(\xi_k-hA_{h,k}(x))+h^2V_{h}(x).
\end{align}
The following proposition, which was proved by \cite{Robert}, provides the existence of the phase function of $h$-FIO's.
%proposition
\begin{proposition}
\label{proposition_Eikonal_1}
Let us fix an open interval $I\Subset(0,\infty)$, $-1<\sigma<1$ and $L>0$.   Then, there exist $R_0,h_0>0$ and a family of smooth and real-valued functions 
$$
\{S_{h}^+;0<h\le h_0,R\ge R_0\}\subset C^\infty(\R^{2d};\R)
$$ 
satisfying the Eikonal equation associated to $\wtilde{p}_{h}$:
\begin{align}
\label{proposition_Eikonal_1_1}
\wtilde{p}_{h}(x,\partial_x S_{h}^+(x,\xi))=|\xi|^2/2,\quad (x,\xi)\in \Gamma^+(R,I,\sigma),
\end{align}
such that
\begin{align}
\label{proposition_Eikonal_1_2_1}
|S_{h}^+(x,\xi)-x\cdot\xi| \le C\<x\>^{1-\mu}\quad,\ x,\xi\in\R^d.
\end{align}
Moreover, for any $|\alpha+\beta|\ge1$,
\begin{align}
\label{proposition_Eikonal_1_2}
|\dderiv{x}{\xi}{\alpha}{\beta}(S_{h}^+(x,\xi)-x\cdot\xi)| &\le C_{\alpha\beta}\min\{R^{1-\mu-|\alpha|},\<x\>^{1-\mu-|\alpha|}\},\quad x,\xi\in\R^d.
\end{align}
Here $C,C_{\alpha\beta}>0$ are independent of $x,\xi,R$ and $h$. 
\end{proposition}

\begin{proof}
Since $h^2V_h$ and $hA_h$ are of long-range type uniformly with respect to $h\in(0,1
]$ (the constant $C_{L,\alpha}$ in \eqref{estimates_on_V_h} can be taken independently of $h$), the proof is same as that of \cite[Proposition 4.1]{Robert} and we omit it. For the $R$ dependence, we refer to \cite[Proposition 3.1]{Bouclet_Tzvetkov_1}.\end{proof}

\begin{remark}
\label{remark_IK}
The crucial point to obtain the estimates \eqref{proposition_Eikonal_1_2_1} and \eqref{proposition_Eikonal_1_2} is the uniform bound \eqref{estimates_on_V_h} and we do not have to use the support properties of $A_h$ and $V_h$. Suppose that  $A$ and $V$ satisfy Assumption \ref{assumption_A} with $\mu\ge0$, \emph{i.e.}, 
$
\<x\>^{-1}|\partial_x^\alpha A(x)|+\<x\>^{-2}|\partial_x^\alpha V(x)|
\le C_{\alpha\beta}\<x\>^{-|\alpha|}
$. 
If $g^{jk} $ is of long-range type, then we still can construct the solution $S_h^+$ to \eqref{proposition_Eikonal_1_1} by using the support properties of $A_h$ and $V_h$, provided that if $L>0$, being independent of $h$, is small enough. However, in this case, $S_h^+-x\cdot\xi$ behaves like $\<x\>^{1-\mu}h^{-1}$ as $h\to 0$, and we cannot obtain the uniform $L^2$-boundedness of the corresponding $h$-FIO. This is one of the reason why we exclude the critical case $\mu=0$. 
\end{remark}

To the phase $S^+_h$ and an amplitude $a\in S(1,g)$, we associate the $h$-FIO defined by
$$
J_h(S_h^+,a)f(x)=(2\pi h)^{-d}\int e^{i(S_h^+(x,\xi)-y\cdot\xi)/h}a(x,\xi)f(y)dyd\xi.
$$
Using \eqref{proposition_Eikonal_1_2}, for sufficiently large $R>0$, we have
\begin{align*}
|\partial_\xi\otimes\partial_x S_h^+(x,\xi)-\Id|&\le C\<R\>^{-\mu}<1/2,\\
|\dderiv{x}{\xi}{\alpha}{\beta}S_h^+(x,\xi)|&\le C_{\alpha\beta}\ \text{for}\ |\alpha+\beta|\ge2,
\end{align*}
uniformly in $h\in(0,h_0]$. Therefore, the standard $L^2$-boundedness of FIO implies that $J_h(S_h^+,a)$ is uniformly bounded on $L^2(\R^d)$ with respect to $h\in(0,h_0]$. 

We now construct the outgoing approximation for $e^{-it\wtilde{H}^h/h}$. 
%theorem
\begin{theorem}	
\label{theorem_IK_1}
Let us fix arbitrarily open intervals $I\Subset I_0\Subset I_1\Subset I_2\Subset(0,\infty)$, $-1<\sigma<\sigma_0<\sigma_1<\sigma_2<1$ and $L>0$. Let $R_0$ and $h_0$ be as in Proposition \ref{proposition_Eikonal_1} with $I,\sigma$ replaced by $I_2,\sigma_2$, respectively. 
Then, for every integer $N\ge0$, the followings hold uniformly with respect to $R\ge R_0$ and $h\in(0,h_0]$.\\
\emph{(1)} There exists a symbol
$$
b^+_h=\sum_{j=0}^{N-1} h^jb^+_{h,j}\quad \text{with}\quad b^+_{h,j}\in S(\<x\>^{-j}\<\xi\>^{-j},g),\ \supp b^+_{h,j}\subset \Gamma^+(R^{1/3},I_1,\sigma_1),
$$
such that, for any $a^+\in S(1,g)$ with $\supp a^+ \subset \Gamma^+(R,I,\sigma)$, we can find
$$
c^+_h=\sum_{j=0}^{N-1} h^jc^+_{h,j}\quad \text{with}\quad c^+_{h,j}\in S(\<x\>^{-j}\<\xi\>^{-j},g),\ \supp c^+_{h,j}\subset \Gamma^+(R^{1/2},I_0,\sigma_0),
$$
such that, for all $0\le t\le h^{-1}$, $e^{-it\wtilde{H}^h/h}\Op_h(a^+)^*$ can be brought to the form
$$
e^{-it\wtilde{H}^h/h}\Op_h(a^+)^*=J_h(S^+_h,b^+_h)e^{-itH_0^h/h}J_h(S^+_h,c^+_h)^*+Q^+_{{\IK}}(t,h,N),
$$
where $J_h(S^+_h,w)$, $w=b^+_h,c^+_h$, are $h$-FIO's associated to the phase $S^+_h$ defined in Proposition \ref{proposition_Eikonal_1} with $R,I$ and $\sigma$ replaced by $R^{1/4},I_2,\sigma_2$, respectively. Moreover, for any integer $s\ge0$ with $2s \le N-1$, 
the remainder $Q^+_{{\IK}}(t,h,N)$ satisfies
\begin{align}
\label{theorem_IK_1_1}
\norm{\<D\>^sQ^+_{{\IK}}(t,h,N)\<D\>^s}_{L^2 \to L^2}\le C_{Ns}h^{N-2s-1},
\end{align}
uniformly with respect to $h\in(0,h_0]$ and $0\le t\le h^{-1}$. \\
\emph{(2)} Let $K_{S^+_h}(t,x,y)$ be the distribution kernel of $J_h(S^+_h,b^+_h)e^{-itH_0^h/h}J_h(S^+_h,c^+_h)^*$. Then, $K_{S^+_h}$ satisfies dispersive estimates:
\begin{align}
\label{theorem_IK_1_2}
|K_{S^+_h}(t,x,y)| \le C|th|^{-d/2},
\end{align}
uniformly with respect to $h\in(0,h_0]$, $x,y\in\R^d$ and $0\le t\le h^{-1}$.
\end{theorem}

\begin{proof} This theorem is basically known and we hence omit the proof. For the construction of the amplitudes $b^+_h$ and $c^+_h$, we refer to Robert \cite[Section 4]{Robert} (see also Bouclet-Tzvetkov \cite[Section 3]{Bouclet_Tzvetkov_1}). The remainder estimate \eqref{theorem_IK_1_1} can be proved by the same argument as that in \cite[Proposition 3.3, Lemma 3.4]{Bouclet_Tzvetkov_1} combined with the following simple estimates:
$$
\norm{\<D\>^s(\wtilde{H}^h+C_1)^{-s/2}}_{L^2 \to L^2}\le C_s h^{-s},\quad s\ge0.
$$
where $C_1>0$ is a large constant. 
Note that these estimates follow from the obvious bounds
$$
\norm{\<D\>^s\<hD\>^{-s}}_{L^2 \to L^2}\le C_sh^{-s},\quad s\ge0,
$$
and the fact that $(\wtilde{p}_h+h\wtilde{p}_{1,h}+C_1)^{-s/2}\in S(\<\xi\>^{-s},g)$ since $\wtilde{p}_h+h\wtilde{p}_{1,h}+C_1$ is uniformly elliptic for sufficiently large $C_1>0$. The dispersive estimate \eqref{theorem_IK_1_2} can be verified by the same argument as that in \cite[Lemma 4.4]{Bouclet_Tzvetkov_1}. 
\end{proof}

The following lemma, which has been essentially proved by \cite{Mizutani2}, tells us that one can still construct the semiclassical approximation for the original propagator $e^{-itH^h/h}$ if we restrict the support of initial data in the region $\Gamma^+(R,J,\sigma)\cap\{x;|x| <h^{-1}\}$. 

%{lemma}
\begin{lemma}
\label{3_lemma_1}
Suppose that $\{a_h^+\}_{h \in (0,1]}$ is a bounded set in $S(1,g)$ with symbols supported in $\Gamma^+(R,I,\sigma)\cap \{x;|x|<h^{-1}\}$. Then, there exists $L>1$ such that, for any $M,s \ge 0$, $h \in (0,h_0]$ and $0 \le t \le h^{-1}$, we have 
$$
\norm{(e^{-itH^h/h}-e^{-it\wtilde{H}^h/h})\Op_h(a_h^+)^*\<D\>^s}_{L^2 \to L^2} \le C_{M,s} h^{M-s},
$$
where $C_{M,s}>0$ is independent of $h$ and $t$.
\end{lemma}

In order to prove this lemma, we need the following. 

\begin{lemma}
\label{lemma_IK_1_1}
Let $f_h\in C^\infty(\R^d)$ be such that for any $\alpha \in \Z^d_+$, 
$$
|\partial_x^\alpha f_h(x)|\le C_\alpha
$$ 
uniformly with respect to $h\in(0,h_0]$ and such that $\supp f_h \subset\{|x|\ge L/(2h)\}$. Let $L>1$ be large enough. Then, under the conditions in Lemma \ref{3_lemma_1}, we have
$$
\norm{f_h(x)\<D\>^\gamma e^{-it\wtilde{H}^h/h}\Op_h(a_h^+)^*\<D\>^s}_{L^2 \to L^2} \le C_{M,s,\gamma} h^{M-s-\gamma},
$$
for any $s,\gamma\ge0$ and $M\ge0$, uniformly with respect $h \in (0,h_0]$ and $0 \le t \le 1/h$.
\end{lemma}

\begin{proof}
We apply Theorem \ref{theorem_IK_1} to $e^{-it\wtilde{H}^h/h}\Op_h(a_h^+)^*$ and obtain
\begin{align*}
e^{-it\wtilde{H}^h/h}\Op_h(a_h^+)^*=J_h(S^+_h,b^+_h)e^{-itH_0^h/h}J_h(S^+_h,c^+_h)^*+Q_{{\IK}}^+(t,h,N).
\end{align*}
By virtue of \eqref{theorem_IK_1_1}, the remainder $f_h(x)\<D\>^\gamma Q_{{\IK}}^+(t,h,N)\<D\>^s$ is bounded on $L^2(\R^d)$ with the norm dominated by $C_{Ns\gamma}h^{N-\gamma-s-1}$, uniformly with respect $h\in(0,h_0]$ and $t\in[0,1/h]$. 
On the other hand, by virtue of \eqref{proposition_Eikonal_1_2}, the phase of $K_{S^+_h}(t,x,y)$, which is given by 
$$
\Phi^+_h(t,x,y,\xi)=S^+_h(x,\xi)-\frac12t|\xi|^2-S^+_h(y,\xi),
$$
satisfies
$
\partial_\xi\Phi^+_h(t,x,y,\xi)=(x-y)(\Id+O(R^{-\mu/4}))-t\xi.
$
We here recall that 
$$
\supp c^+_h\subset \{(y,\xi)\in \R^{2d};a_h^+(y,\partial_\xi S^+_h(y,\xi))\neq 0\}
$$ 
(see, \cite[Lemma 3.2]{Mizutani2} and its proof). In particular, $c^+_h(y,\xi)$ vanishes in the region $\{y;|y| \ge 1/h\}$. 
We now set $L=4\sqrt{\sup I_2}+2$, where $I_2$ is given in Theorem \ref{theorem_IK_1}. Since $|x|\ge L/(2h)$, $ |y|<1/h$ and $|\xi|^2 \in I_2$ 
on the support of the amplitude $f_h(x)b^+_h(x,\xi)\overline{c^+_h(y,\xi)}$, we obtain
$$
|\partial_\xi\Phi^+_h(t,x,y,\xi)| >c(1+|x|+|y|+|\xi|+t+h^{-1}),\quad 0 \le t \le h^{-1},
$$
for some universal constant $c>0$. The assertion now follows from an integration by parts and  the $L^2$-boundedness of $h$-FIO's. 
\end{proof}

\begin{proof}[Proof of Lemma \ref{3_lemma_1}] 
The Duhamel formula yields
\begin{align*}
(e^{-itH^h/h}-e^{-it\wtilde{H}^h/h})
&=-\frac ih\int_0^t e^{-i(t-s)H^h/h}W_0^he^{-is\wtilde{H}^h/h}ds\\
&=-\frac ih\int_0^t e^{-i(t-s)H^h/h}e^{-is\wtilde{H}^h/h}W_0^hds\\
&\ \ \ +\frac{1}{h^2}\int_0^t e^{-i(t-s)H^h/h}\int_0^s e^{-i(s-\tau)\wtilde{H}^h/h}[\wtilde{H}^h,W_0^h]e^{-i\tau \wtilde{H}^h/h}d\tau ds,
\end{align*}
where $W_0^h:=H^h-\wtilde{H}^h$ consists of the following two parts:
\begin{align*}
&\frac{i h^2}{2}\sum_{j,k}\left(\partial_j g^{jk}(1-\psi(hx/L))A_k+ (1-\psi(hx/L))A_jg^{jk}\partial_k\right),\\
&\frac{h^2}{2}\sum_{j,k}(1-\psi(hx/L))^2g^{jk}A_jA_k+h^2(1-\psi(hx/L))V.
\end{align*}
In particular, $W_0^h$ is a first order differential operator of the form
$$
h^2\sum_{|\alpha|=1}f_{\alpha}^h(x)\partial_x^\alpha+h^2f_0^h(x),
$$
where $f_\alpha^h,f_0^h$ are supported in $\{|x|\ge L/(2h)\}$ and satisfy
\begin{align}
\label{proof_lemma_1_1}
|\partial_x^\beta f^h_\alpha(x)|\le C_{\alpha\beta}\<x\>^{1-\mu-|\beta|},\ 
|\partial_x^\beta f_0^h(x)|\le C_{\alpha\beta}\<x\>^{2-\mu-|\beta|}.
\end{align}
Since $\{|x|\ge L/(2h)\}\cap\pi_x(\supp a_h^+)=\emptyset$ if $L>1$, we have
$$
\norm{W_0^h\Op_h(a_h^+)^*\<D\>^s}_{L^2 \to L^2} \le C_{M,s}h^{M-s},\quad M\ge0,\ s\in\R.
$$
Therefore, the first term of the right hand side of the above Duhamel formula satisfies the desired estimates since $e^{-itH^h/h}$ and $e^{-it\what{H}^h/h}$ are unitary on $L^2$. 

We next study the second term. 
Again by the Duhamel formula, we have
$$
[\wtilde{H}^h,W_0^h]e^{-i\tau \wtilde{H}^h/h}
=e^{-i\tau \wtilde{H}^h/h}[\wtilde{H}^h,W_0^h]+\frac ih\int_0^\tau e^{-i(\tau-u)\wtilde{H}^h/h}[\wtilde{H}^h,[\wtilde{H}^h,W_0^h]]e^{-iu \wtilde{H}^h/h}du.
$$
Since the coefficients of the commutator $[\wtilde{H}^h,W_0^h]$ are supported in $\{|x|\ge L/(2h)\}$, the support property of $a^+_h$ again implies that $[\wtilde{H}^h,W_0^h]\Op_h(a_h^+)^*\<D\>^s=O_{L^2\to L^2}(h^{M-s})$ for any $M\ge0$ and $s\in\R$. Furthermore, by virtue of \eqref{estimates_on_V_h}, \eqref{proof_lemma_1_1} and the symbolic calculus, the coefficients of $[\wtilde{H}^h,[\wtilde{H}^h,W_0^h]]$ are uniformly bounded in $x$ and supported in $\{|x|\ge L/(2h)\}$. We now apply Lemma \ref{lemma_IK_1_1} to $[\wtilde{H}^h,[\wtilde{H}^h,W_0^h]]e^{-iu\wtilde{H}^h/h}\Op_h(a_h^+)^*$ and obtain the assertion.
\end{proof}

We now come into the proof of Proposition \ref{theorem_dispersive_IK}. 
\begin{proof}[Proof of Proposition \ref{theorem_dispersive_IK}]
Rescaling $t\to th$, it suffices to show
$$
\norm{\Op_h(a_h^+)e^{-itH^h/h}\Op_h(a_h^+)^*}_{L^1\to L^\infty} \le C_{\ep}|th|^{-d/2},\quad 0<  |t| \le h^{-1},
$$
where $H^h=h^2H$. Let $A_h(x,y)$ be the distribution kernel of $\Op_h(a_h^+)$:
$$A_h(x,y)=(2\pi h)^{-d}\int e^{i(x-y)\cdot\xi/h}a_h^+(x,\xi)d\xi.$$
Since $a_h^+\in S(1,g)$ is compactly supported in $I$ with respect to $\xi$, we easily see that
$$
\sup_x\int |A_h(x,y)|dy+\sup_y\int |A_h(x,y)|dx\le C,\quad h\in(0,1].
$$
Moreover, since $\<\xi\>^sa_h^+\<\xi\>^{\gamma}\in S(1,g)$ for any $s,\gamma$, we have 
\begin{align}
\label{6}
\norm{\<D\>^s\Op_h(a_h^+)\<D\>^{\gamma}}_{L^2 \to L^2}
\le C_sh^{-s-\gamma}.
\end{align}
Combining these two estimates with Theorem \ref{theorem_IK_1} and Lemma \ref{3_lemma_1}, we can write
$$
\Op_h(a_h^+)e^{-itH^h/h}\Op_h(a_h^+)^*=K_1(t,h,N)+K_2(t,h,N),
$$
where 
\begin{align*}
K_1(t,h,N)&=\Op_h(a_h^+)J_h(S^+_h,b^+_h)e^{-itH_0^h/h}J_h(S^+_h,c^+_h)^*,\\
K_2(t,h,N)&=\Op_h(a_h^+)Q^+_{{\IK}}(t,h,N)+\Op_h(a_h^+)(e^{-itH^h/h}-e^{-it\wtilde{H}^h/h})\Op_h(a_h^+)^*.
\end{align*} 
By \eqref{theorem_IK_1_2}, the distribution kernel of $K_1(t,h,N)$, which we denote by $K_1(t,x,y)$, satisfies
$$
|K_1(t,x,y)| \le \int |A_h(x,z)||K_{S^+_h}(t,z,y)|dz \le C_N|th|^{-d/2},\quad 0<t\le h^{-1},
$$
uniformly in $h \in(0,h_0]$. On the other hand, \eqref{theorem_IK_1_1}, Lemma \ref{3_lemma_1} and \eqref{6} imply
$$
\norm{\<D\>^s K_2(t,h,N)\<D\>^s}_{L^2 \to L^2} \le C_{N,s} h^{N-2s-1}.
$$
If we choose $N\ge d+2$ and $s>d/2$, then it follows from the Sobolev embedding that the distribution kernel of  $K_2(t,h,N)$ is uniformly bounded in $\R^{2d}$ with respect to $h\in(0,h_0]$ and $0<t\le h^{-1}$. Therefore, $\Op_h(a_h^+)e^{-itH^h/h}\Op_h(a_h^+)^*$ has the distribution kernel $K(t,x,y)$ satisfying dispersive estimates for $0<t\le h^{-1}$:
\begin{align}
\label{5}
|K(t,x,y)| \le C_N|th|^{-d/2},\quad x,y\in\R^d.
\end{align}
Finally, using the following relation,
$$
\Op_h(a_h^+)e^{-itH^h/h}\Op_h(a_h^+)^*=\left(\Op_h(a_h^+)e^{itH^h/h}\Op_h(a_h^+)^*\right)^*,
$$
we learn $K(t,x,y)=\overline{K(-t,y,x)}$ and \eqref{5} also holds for $0<-t\le h^{-1}$. For the incoming case, the proof is analogous and we omit it. 
\end{proof}

%%%%%%%%%%%%%%%%%%%%%%%%%%%%%%%%%%%%%%%%%%%%%%%%      Parametrix construction in terms of time dependent Fourier integral operators	%%%%%%%%%%%%%%%%%%%%%%%%%%%%%%%%%%%%%%%%%%%%%%

\section{Fourier integral operators with the time dependent phase}
\label{section_WKB}
Throughout this section we assume Assumption \ref{assumption_A} with $\mu\ge0$. 
Consider a symbol $\chi_\ep \in S(1,g)$ supported in a region
$$
\Omega(\ep):=\{(x,\xi)\in\R^{2d}; \<x\>>\ep|\xi|/2\},
$$
where $\ep>0$ is an arbitrarily small fixed constant. In this section we prove the following dispersive estimate:
$$
\norm{\Op(\chi_\ep)e^{-itH}\Op(\chi_\ep)^*}_{L^1\to L^\infty}\le C_\ep |t|^{-d/2},\quad 0<|t|\le t_\ep,
$$
where $t_\ep>0$ is a small constant depending on $\ep$. This estimate, combined with the $L^2$-boundedness of $\Op(\chi_\ep)$ and $e^{-itH}$, implies Strichartz estimates for $\Op(\chi_\ep)e^{-itH}$. 

Let us give a short summary of the steps of proof. Choose $\chi_\ep^*\in S(1,g)$ so that $\supp \chi_\ep^*=\supp\chi_\ep$ and $\Op(\chi_\ep)^*=\Op(\chi_\ep^*)+\Op(r_N)$ with some $r_N\in S(\<x\>^{-N}\<\xi\>^{-N},g)$ for sufficiently large $N>d/2$. We first construct an approximation for $e^{-itH}\Op(\chi_\ep^*)$ in terms of the FIO with a time dependent phase:
$$
J(\Psi,b^N)f(x)=\frac{1}{(2\pi)^d}\int e^{i(\Psi(t,x,\xi)-y\cdot\xi)}b(t,x,\xi)f(y)dyd\xi,
$$
where $\Psi$ is a generating function of the Hamilton flow associated to $p(x,\xi)$ and $(\partial_\xi \Psi,\xi)\mapsto (x,\partial_x\Psi)$ is the corresponding canonical map, and the amplitude $b= b_0+b_2+\cdots+b_{N-1}$ solves the corresponding transport equations. Although such parametrix constructions are well known as WKB approximations (at least if $\chi_\ep^*$ is compactly supported in $\xi$ and the time scale depends on the size of frequency), we give the detail of proof since, in the present case, $\supp \chi_\ep^*$ is not compact with respect to $\xi$ and $t_\ep$ is independent of the size of frequency. The crucial point is that $p(x,\xi)$ is of \emph{quadratic type} on $\Omega(\ep)$:
$$
|\dderiv{x}{\xi}{\alpha}{\beta}p(x,\xi)|\le C_{\alpha\beta},\quad (x,\xi)\in\Omega(\ep),\ |\alpha+\beta|\ge2,
$$
which allows us to follow a classical argument (due to, \emph{e.g.}, \cite{Kitada_Kumanogo}) and construct the approximation for $|t|<t_\ep$ if $t_\ep>0$ is small enough. The composition $\Op(\chi_\ep)J(\Psi,b)$ is also a FIO with the same phase, and a standard stationary phase method can be used to prove dispersive estimates for $0<|t|<t_\ep$. It remains to obtain the $L^1 \to L^\infty$ bounds of the remainders $\Op(\chi_\ep)e^{-itH}\Op(r_N)$ and $\Op(\chi_\ep)e^{-itH}(\Op(\chi_\ep^*)-J(\Psi,b^N))$. If $e^{-itH}$ maps from the Sobolev space $H^{d/2}(\R^d)$ to itself, then $L^1 \to L^\infty$ bounds are direct consequences of the Sobolev embedding and $L^2$-boundedness of PDO. However, our Hamiltonian $H$ is not bounded below (on $\{|x|\gtrsim|\xi|\}$) and such a property does not hold in general. To overcome this difficulty, we use an Egorov type lemma as follows. By the Sobolev embedding and the Littlewood-Paley decomposition, the proof is reduced to that of the following estimates:
\begin{align}
\label{sketch_WKB}
\sum_{j\ge0}\norm{2^{j\gamma}S_j(D)\Op(\chi_\ep)e^{-itH}\Op(r_N)\<D\>^{\gamma}f}_{L^2}^2\le C \norm{f}_{L^2}^2,
\end{align}
where $\gamma>d/2$ and $S_j$ is a dyadic partition of unity. Then, we will prove that there exists $\eta_j(t,\cdot,\cdot)\in S(1,g)$ such that
$$
2^j\le C (1+|x|+|\xi|)\quad \text{on}\ \supp\eta_j(t),
$$
and that
$$
S_j(D)\Op(\chi_\ep)e^{-itH}=e^{-itH}\Op(\eta_j(t))+O_{L^2\to L^2}(2^{-jN}),\quad |t|<t_\ep\ll1.
$$  
Choosing $\delta>0$ with $\gamma+\delta\le N/2$, we learn that $2^{j(\gamma+\delta)}\eta_j(t)r_N\<\xi\>^{\gamma} \in S(1,g)$ and hence \eqref{sketch_WKB}. $\Op(\chi_\ep)e^{-itH}(\Op(\chi_\ep^*)-J(\Psi,b))$ can be controlled similarly. 

%Hamilton flow generated by the total energy

\subsection{Short-time behavior of the Hamilton flow}
\label{subsection_WKB_1}
This subsection discusses the classical mechanics generated by $p(x,\xi)$. We denote the solution to the following Hamilton equations by $(X(t),\Xi(t))=(X(t,x,\xi),\Xi(t,x,\xi))$:
\begin{equation}
\left\{
\begin{aligned}
\nonumber
\dot{X}_j&=\frac{\partial p}{\partial \xi_j}(X,\Xi)=\sum_k g^{jk}(X)(\Xi_k-A_k(X)),\\
\dot{\Xi}_j&=-\frac{\partial p}{\partial x_j}(X,\Xi)
=-\frac12 \sum_{k,l}\frac{\partial g^{kl}}{\partial x_j}(X)(\Xi_k-A_k(X))(\Xi_l-A_l(X))\\
&\quad\quad\quad\quad\quad\quad\quad\ \ +\sum_{k,l}g^{kl}(X)\frac{\partial A_k}{\partial x_j}(X)(\Xi_l-A_l(X))
-\frac{\partial V}{\partial x_j}(X)
\end{aligned}
\right.
\end{equation}
with the initial condition $(X(0),\Xi(0))=(x,\xi)$, where $\dot{f}=\partial_t f$. We first observe that the flow conserves the energy:
$$p(x,\xi)=p(X(t),\Xi(t)),$$
which, combined with the uniform ellipticity of $g^{jk}$, implies 
\begin{align*}
|\Xi(t)-A(X(t))|^2
&\lesssim p(X(t),\Xi(t))-V(X(t))\\
&=p(x,\xi)-V(X(t))\\
&\lesssim |\xi-A(x)|^2+|V(x)|+|V(X(t))|,
\end{align*}
and hence $|\Xi(t)|\lesssim|\xi|+\<x\>+\<X(t)\>$. By the Hamilton equation, we then have
\begin{align*}
|\dot{X}(t)|+|\dot{\Xi}(t)|\le C (1+|\xi|+|x|+|X(t)|+|\Xi(t)|).
\end{align*}
Applying Gronwall's inequality to this estimate, we obtain an a priori bound:
$$
|X(t)-x|+|\Xi(t)-\xi| \le C_T |t|(1+|x|+|\xi|),\quad |t|\le T,\ x,\xi\in\R^d. 
$$
Using this estimate, we obtain more precise behavior of the flow with initial conditions in $\Omega(\ep)$.

%\lemma
\begin{lemma}
\label{WKB_1_lemma_3}
Let $\ep>0$. Then, for sufficiently small $t_\ep>0$ and all $\alpha,\beta\in \Z^d_+$,
$$
|\dderiv{x}{\xi}{\alpha}{\beta}(X(t,x,\xi)-x)| +|\dderiv{x}{\xi}{\alpha}{\beta}(\Xi(t,x,\xi)-\xi|
\le C_{\alpha\beta\ep}|t| \<x\>^{1-|\alpha+\beta|},
$$
uniformly with respect to $(t,x,\xi)\in(-t_\ep,t_\ep)\times\Omega(\ep)$.
\end{lemma}

\begin{proof}
We only consider the case with $t \ge 0$, the proof for the opposite case is similar. Let $(x,\xi)\in\Omega(\ep)$. At first we remark that for sufficiently small $t_\ep>0$, 
\begin{align}
\label{WKB_1_lemma_1_1}
\<x\>/2 \le |X(t,x,\xi)| \le 2\<x\>,\quad |t|\le t_\ep.
\end{align}
For $|\alpha+\beta|=0$, the assertion is obvious. We let $|\alpha+\beta|=1$ and differentiate the Hamilton equations with respect to $\partial_x^\alpha\partial_\xi^\beta$: 
\begin{align}
\label{WKB_1_lemma_3_3}%\label{WKB_1_lemma_3_3}
\frac{d}{dt}
\left(\begin{matrix}
\dderiv{x}{\xi}{\alpha}{\beta} X\\
\dderiv{x}{\xi}{\alpha}{\beta} \Xi
\end{matrix}\right)
=
\left(\begin{matrix}
\partial_{x}\partial_\xi p(X,\Xi) & \partial_{\xi}^2p(X,\Xi)\\
-\partial_{x}^2p(X,\Xi) & -\partial_{\xi}\partial_{x}p(X,\Xi)
\end{matrix}\right)
\left(\begin{matrix}
\dderiv{x}{\xi}{\alpha}{\beta} X\\
\dderiv{x}{\xi}{\alpha}{\beta} \Xi
\end{matrix}\right).
\end{align} 
Using \eqref{WKB_1_lemma_1_1}, we learn that $p(X(t),\Xi(t))$ is of quadratic type in $\Omega(\ep)$:
\begin{align*}
\abs{(\dderiv{x}{\xi}{\alpha}{\beta}p)(X(t),\Xi(t))}\le
C_{\alpha\beta\ep}\<x\>^{2-|\alpha+\beta|},\quad (t,x,\xi)\in(-t_\ep,t_\ep)\times\Omega(\ep).
\end{align*}
All entries of the above matrix hence are uniformly bounded in $(t,x,\xi)\in(-t_\ep,t_\ep)\times\Omega(\ep)$. Taking $t_\ep>0$ smaller if necessary, integrating \eqref{WKB_1_lemma_3_3} with respect to $t$ and applying Gronwall's inequality, we have the assertion with $|\alpha+\beta|=1$. For $|\alpha+\beta| \ge 2$, we prove the estimate for $\partial_{\xi_1}^2X(t)$ and $\partial_{\xi_1}^2\Xi(t)$ only, where $\xi=(\xi_1,\xi_2,...,\xi_d)$. Proofs for other cases are similar, and proofs for higher derivatives follow from an induction on $|\alpha+\beta|$. By the Hamilton equation, we learn
\begin{align*}
\frac{d}{dt}\partial_{\xi_1}^2X(t)
=\partial_{x}\partial_\xi p(X(t),\Xi(t)) \partial_{\xi_1}^2X(t)
+\partial_{\xi}^2p(X(t),\Xi(t)) \partial_{\xi_1}^2\Xi(t)+Q(X(t),\Xi(t)),
\end{align*}
where $Q(X(t),\Xi(t))$ satisfies
\begin{align*}
|Q(X(t),\Xi(t))| 
&\le C_\ep\sum_{|\alpha+\beta|=3,|\beta| \ge 1}
|(\partial_x^\alpha\partial_\xi^\beta p)(X(t),\Xi(t))|
|\partial_{\xi_1}X(t)|^{|\alpha|}|\partial_{\xi_1}\Xi(t)|^{|\beta|}\\
&\le C_\ep\<x\>^{-1}. 
\end{align*}
We similarly obtain
\begin{align*}
\frac{d}{dt}\partial_{\xi_1}^2\Xi(t)
=-\partial_{x}^2p(X(t),\Xi(t)) \partial_{\xi_1}^2X(t)
-\partial_{\xi}\partial_{x}p(X(t),\Xi(t)) \partial_{\xi_1}^2\Xi(t)
+O(\<x\>^{-1}).
\end{align*}
Applying Gronwall's inequality, we have the desired estimates.
 \end{proof}

%lemma
\begin{lemma}
\label{WKB_2_lemma_1}
\emph{(1)} Let $t_\ep>0$ be small enough. Then, for any $|t| < t_\ep$, the map
$$
g(t):(x,\xi) \mapsto (X(t,x,\xi),\xi)
$$
is a diffeomorphism from $\Omega(\ep/2)$ onto its range, and satisfies
$$\Omega(\ep)\subset g(t,\Omega(\ep/2))\ \text{for all}\ |t|<t_\ep.$$
\emph{(2)} Let $\Omega(\ep) \ni (x,\xi)\mapsto(Y(t,x,\xi),\xi)\in \Omega(\ep/2)$ be the inverse map of $g(t)$. Then, $Y(t,x,\xi)$ and $\Xi(t,Y(t,x,\xi),\xi)$ satisfy the same estimates as that for $X(t,x,\xi)$ and $\Xi(t,x,\xi)$ of Lemma \ref{WKB_1_lemma_3}, respectively:
$$
|\dderiv{x}{\xi}{\alpha}{\beta}(Y(t,x,\xi)-x)| +|\dderiv{x}{\xi}{\alpha}{\beta}(\Xi(t,Y(t,x,\xi),\xi)-\xi|
\le C_{\alpha\beta\ep}|t| \<x\>^{1-|\alpha+\beta|},
$$
uniformly with respect to $(t,x,\xi)\in(-t_\ep,t_\ep)\times\Omega(\ep)$.
\end{lemma}

\begin{proof}
Choosing a cut-off function $\rho \in S(1,g)$ such that $0 \le \rho \le 1$, $\supp \rho \subset \Omega(\ep/3)$ and $\rho \equiv 1$ on $\Omega(\ep/2)$, we modify $g(t)$ as follows:
$$
g_\rho(t,x,\xi)=(X_\rho(t,x,\xi),\xi),\ X_\rho(t,x,\xi)=(1-\rho(x,\xi))x+\rho(x,\xi)X(t,x,\xi).
$$ 
It is easy to see that, for $(t,x,\xi)\in(-t_\ep,t_\ep)\times\Omega(\ep/2)$, $g_\rho(t,x,\xi)$ is  smooth and Lemma \ref{WKB_1_lemma_3} implies 
\begin{align*}
|\dderiv{x}{\xi}{\alpha}{\beta}g_\rho(t,x,\xi)|&\le C_{\alpha\beta\ep},\quad |\alpha+\beta|\ge 1,\\
|J(g_\rho)(t,x,\xi)-\Id| &\le C_\ep t_\ep,
\end{align*}
where $J(g_\rho)$ is the Jacobi matrix with respect to $(x,\xi)$ and the constant $C_\ep>0$ is independent of $t,x$ and $\xi$. Choosing $t_\ep>0$  so small that $C_\ep t_\ep<1/2$, and applying the Hadamard global inverse mapping theorem, we see that, for any fixed $|t| <t_\ep$, $g_\rho(t)$ is a diffeomorphism from $\R^{2d}$ onto itself. By definition, $g(t)$ is diffeomorphic from $\Omega(\ep/2)$ onto its range. Since $g_\rho(t)$ is bijective, it remains to check that
$$
\Omega(\ep)^c \supset g_\rho(t,\Omega(\ep/2)^c),\quad |t|<t_\ep.
$$
Suppose that $(x,\xi) \in \Omega(\ep/2)^c$. If $(x,\xi)  \in \Omega(\ep/3)^c$, then the assertion is obvious since $g_\rho(t)\equiv \Id$ outside $\Omega(\ep/3)$. If $(x,\xi) \in \Omega(\ep/3) \setminus \Omega(\ep/2)$, then, by Lemma \ref{WKB_1_lemma_3} and the support property of $\rho$, we have
\begin{align*}
|X_\rho(t,x,\xi)|
\le |x|+\rho(x,\xi)|(X(t,x,\xi)-x)|
\le (\ep/2+C_0t_\ep)\<\xi\>
\end{align*}
for some $C_0>0$ independent of $x,\xi$ and $t_\ep$. Choosing $t_\ep<\ep/(2C_0)$, we obtain the assertion.

We next prove the estimates on $Y(t)$. Since $(Y(t,x,\xi),\xi) \in \Omega(\ep/2)$, we learn
\begin{align*}
|Y(t,x,\xi)-x|
&= |X(0,Y(t,x,\xi),\xi)-X(t,Y(t,x,\xi),\xi)|\\
&\le \sup_{(x,\xi)\in\Omega(\ep/2)}|X(t,x,\xi)-x| \\
&\le C_\ep |t| \<x\>. 
\end{align*}
For $\alpha,\beta \in \Z_+^d$ with $|\alpha+\beta|=1$, apply $\dderiv{x}{\xi}{\alpha}{\beta}$ to the equality $x=X(t,Y(t,x,\xi),\xi)$. We then have the following equality
\begin{align*}
A(t,Z(t,x,\xi))
\dderiv{x}{\xi}{\alpha}{\beta}(Y(t,x,\xi)-x)
=\dderiv{y}{\eta}{\alpha}{\beta}(y-X(t,y,\eta))
|_{(y,\eta)=Z(t,x,\xi)},\
\end{align*}
where $Z(t,x,\xi)=(Y(t,x,\xi),\xi)$ and $A(t,Z)=(\partial_x X)(t,Z)$ is a $d\times d$-matrix. By Lemma \ref{WKB_1_lemma_3} and a similar argument as that in the proof of Lemma \ref{WKB_2_lemma_1} (1), we learn that $A(t,Z(t,x,\xi))$ is invertible if $t_\ep>0$ is small enough, and that $A(t,Z(t,x,\xi))$ and $A(t,Z(t,x,\xi))^{-1}$ are bounded uniformly in $(t,x,\xi) \in (-t_\ep,t_\ep)\times\Omega(\ep/2)$ . Therefore, 
\begin{align*}
\abs{\dderiv{x}{\xi}{\alpha}{\beta}(Y(t,x,\xi)-x)} 
&\le C_{\alpha\beta}\sup_{(x,\xi)\in\Omega(\ep/2)}|\dderiv{x}{\xi}{\alpha}{\beta}(x-X(t,x,\xi))|\\
&\le C_{\alpha\beta}|t| \<x\>^{1-|\alpha+\beta|}. 
\end{align*}
Proofs for higher derivatives are obtained by an induction with respect to $|\alpha+\beta|$ and proofs for $\Xi(t,Y(t,x,\xi),\xi)$ are similar. 
\end{proof}

%\subsection{The parametrix for $\Op(\chi_\ep)e^{-itH}\Op(\chi_\ep)^*$}

\subsection{The parametrix for $\Op(\chi_\ep)e^{-itH}\Op(\chi_\ep)^*$}
Before starting the construction of parametrix, we prepare two lemmas. The following is an Egorov type theorem which will be used to control the remainder term. We write $\exp tH_p(x,\xi)=(X(t,x,\xi),\Xi(t,x,\xi))$.

%{lemma}
\begin{lemma}							
\label{lemma_Egorov}
For $h\in(0,1]$, consider a $h$-dependent symbol $\eta_h \in S(1,g)$ such that $\supp \eta_h \subset \Omega(\ep)\cap\{1/(2h)<|\xi|<2/h\}$.  
Then, for sufficiently small $t_\ep>0$, independent of $h$, and any integer $N\ge0$, there exists a bounded family of symbols $\{\eta^{N}_h(t,\cdot,\cdot);|t|<t_\ep,\ 0<h\le1\}\subset S(1,g)$ such that
$$
\supp \eta^{N}_h(t,\cdot,\cdot)\subset \exp (-t)H_p(\supp \eta_h),
$$
and that
$$
\norm{e^{itH}\Op(\eta_h)e^{-itH}-\Op(\eta^{N}_h(t))}_{L^2\to L^2}\le C_{N\ep}h^N,
$$ 
uniformly with respect to $0<h\le1$ and $|t|<t_\ep$.
\end{lemma}

\begin{proof}
Let $\eta_h^0(t,x,\xi)=\eta_h(\exp tH_p(x,\xi))=\eta_h(X(t,x,\xi),\Xi(t,x,\xi))$. It is easy to see that $\supp \eta^0_h \subset \exp (-t)H_p(\supp \eta_h)$. Moreover, Lemma \ref{WKB_1_lemma_3} implies that $\{\eta^0_h; |t|<t_\ep,0<h\le1\}$ is a bounded subset of $S(1,g)$. By a direct computation, $\eta^0_h$ solves 
$$
\partial_t \eta^0_h=\{p,\eta^0_h\};\quad \eta^0_h|_{t=0}=\eta_h,
$$
where $\{\cdot,\cdot\}$ is the Poisson bracket. Then, by a standard pseudodifferential calculus, there exists a bounded set $\{r^0_h(t,\cdot,\cdot); 0\le t<t_\ep,0<h\le1\}\subset S(1,g)$ with $\supp r^0_h \subset \exp(-t)H_p(\supp \eta_h)$ such that
$$
\frac{d}{dt} \Op(\eta^0_h)=i[H,\Op(\eta^0_h)]+h\Op(r^0_h).
$$
We next set 
$$
\eta^1_h(t,x,\xi)=\int_0^tr^0_h(s,X(t-s,x,\xi),\Xi(t-s,x,\xi))ds.
$$
Again, we learn that $\{\eta^1_h(t,\cdot,\cdot);|t|<t_\ep,0<h\le1\} \subset S(1,g)$ is also bounded and that $\supp \eta^1_h \subset \exp(-t)H_p(\supp \eta_h)$ for all $|t|<t_\ep$ and $0<h\le1$. Moreover, $\eta^1_h$ solves
$$
\partial_t \eta^1_h=\{p,\eta^1_h\}+r^0_h;\quad \eta^1_h|_{t=0}=0,
$$
which implies 
$$
\frac{d}{dt} \Op(\eta^0_h+h\eta^1_h)=i[H,\Op(\eta^0_h+h\eta^1_h)]+h^2\Op(r^1_h).
$$
with some $\{r^1_h;0\le t<t_\ep,0<h\le1\}\subset S(1,g)$ and $\supp r^1_h \subset \exp(-t)H_p(\supp \eta_h)$. Iterating this procedure and putting $\eta^{N}_h=\sum_{j=0}^{N-1}h^{j}\eta_h^j$, we obtain the assertion. 
\end{proof}

Using this lemma, we have the following.
%{lemma}
\begin{lemma}
\label{lemma_WKB}Let $\ep>0$. Then, for any symbol $\chi_\ep \in S(1,g)$ with $\supp \chi_\ep \subset \Omega(\ep)$ and any integer $N\ge1$, there exists $\chi_\ep^*\in S(1,g)$ with $\supp \chi_\ep^*\subset \Omega(\ep)$ such that for any $\gamma<N/2$,
$$
\sup_{|t|<t_\ep}\norm{\Op(\chi_\ep)e^{-itH}\Op(\chi_\ep)^*-\Op(\chi_\ep)e^{-itH}\Op(\chi_\ep^*)}_{H^{-\gamma}(\R^d)\to H^{\gamma}(\R^d)}
\le C_{N\gamma\ep}
$$
\end{lemma}

\begin{proof}
By the expansion formula \eqref{asymptotic_expansion_adjoint}, there exists $\chi_\ep^*\in S(1,g)$ with $\supp \chi_\ep^*\subset \Omega(\ep)$ such that 
$$\Op(\chi_\ep)^*=\Op(\chi_\ep^*)+\Op(r_0(N))$$ 
with some $r_0(N)\in S(\<x\>^{-N}\<\xi\>^{-N},g)$. For $\delta>0$ with $2\gamma+\delta\le N$, we split 
\begin{align*}
&\<D\>^\gamma\Op(\chi_\ep)e^{-itH}\Op(r_0(N))\<D\>^{\gamma}\\
&=\<D\>^\gamma\Op(\chi_\ep)e^{-itH}\<D\>^{-\gamma-\delta}\<x\>^{-\gamma-\delta}\cdot\<x\>^{\gamma+\delta}\<D\>^{\gamma+\delta}\Op(r_0(N))\<D\>^{\gamma}.
\end{align*}
Since $\<x\>^{\gamma+\delta}\<\xi\>^{\gamma+\delta} r_0(N)\<\xi\>^\gamma \in S(1,g)$, $\<x\>^{\gamma+\delta}\<D\>^{\gamma+\delta}\Op(r_0(N))\<D\>^{\gamma}$ is bounded on $L^2$. 
In order to prove the $L^2$-boundedness of the first term of the right hand side, we use the standard Littlewood-Paley decomposition and Lemma \ref{lemma_Egorov} as follows. 
Consider a dyadic partition of unity with respect to the frequency:
$$\sum_{j=0}^\infty S_j(D)=1,$$
where $S_j(\xi)=S(2^{-j}\xi)$, $j\ge1$, with some $S\in C_0^\infty(\R^d)$ supported in $\{1/2<|\xi|<2\}$ and $S_0\in C^\infty_0(\R^d)$ supported in $\{|\xi|<1\}$. Then, 
\begin{align*}
&\norm{\<D\>^\gamma\Op(\chi_\ep)e^{-itH}\<D\>^{-\gamma-\delta}\<x\>^{-\gamma-\delta}f}_{L^2}\\
&\le C\bigg(\sum_{j=0}^\infty\norm{2^{j\gamma}S_j(D)\Op(\chi_\ep)e^{-itH}\<D\>^{-\gamma-\delta}\<x\>^{-\gamma-\delta}f}_{L^2}^2\bigg)^{1/2}.
\end{align*}
By the expansion formula \eqref{asymptotic_expansion_composition}, there exists a sequence of symbols $\eta_j\in S(1,g)$ supported in $\Omega(\ep)\cap\{2^{j-1}<|\xi|<2^{j+1}\}
$ such that
$$
S_j(D)\Op(\chi_\ep)=\Op(\eta_j)+Q_1(j,N),\quad\norm{Q_1(j,N)}_{L^2\to L^2}=O(2^{-jN}).
$$
We then learn by Lemma \ref{lemma_Egorov} with $h=2^{-j}$ that there exists $\{\eta^N_j(t);|t|<t_\ep\}\subset S(1,g)$ such that
$$
\Op(\eta_j)e^{-itH}=e^{-itH}\Op(\eta^N_j(t))+Q_2(t,j,N),\quad \sup_{|t|<t_\ep}\norm{Q_2(t,j,N)}_{L^2\to L^2}=O(2^{-jN}).
$$
Since $N \ge \gamma+\delta$, the remainder satisfies
$$
\sup_{|t|<t_\ep}\norm{2^{j\gamma}\left(Q_1(j,N)e^{-itH}+Q_2(t,j,N)\right)\<D\>^{-\gamma-\delta}\<x\>^{-\gamma-\delta}f}_{L^2}^2\le C2^{-2j\delta}\norm{f}_{L^2}^2.
$$
Suppose that $(x,\xi)\in \supp \eta^N_j(t)$. Since $\supp \eta^N_j(t) \subset \exp(-t)H_p(\supp \eta_j)$, we have
$$|X(t,x,\xi)|>\ep\<\Xi(t,x,\xi)\>,\ 2^{j-1}<|\Xi(t,x,\xi)|<2^{j+1}.$$
Using Lemma \ref{WKB_1_lemma_3} with the initial data $(X(t,x,\xi),\Xi(t,x,\xi))$, we learn
$$
|x-X(t,x,\xi)|+|\xi-\Xi(t,x,\xi)| \le Ct_\ep\<X(t,x,\xi)\>, \quad |t|<t_\ep.
$$
Combining these two estimates, we see that
$$
2^j\le C( 1+|x|+|\xi|),\quad (x,\xi)\in \supp\eta^N_j(t),\ |t|<t_\ep,
$$
where the constant $C>0$ is independent of $x$, $\xi$ and $t$, provided that $t_\ep>0$ is small enough. Therefore, 
$
2^{j(\gamma+\delta)}\eta^N_j(t)\<\xi\>^{-\gamma-\delta}\<x\>^{-\gamma-\delta}\in S(1,g)
$  
and the corresponding PDO is bounded on $L^2$. Finally, we obtain
\begin{align*}
&\sum_{j=0}^\infty\norm{2^{j\gamma}\Op(\eta_j)e^{-itH}\<D\>^{-\gamma-\delta}\<x\>^{-\gamma-\delta}f}_{L^2}^2\\
&\le C\sum_{j=0}^\infty\left(\norm{2^{-j\delta}2^{j(\gamma+\delta)}\Op(\eta^N_j(t))\<D\>^{-\gamma-\delta}\<x\>^{-\gamma-\delta}f}_{L^2}^2+2^{-2j\delta}\norm{f}_{L^2}^2\right)\\
&\le C\sum_{j=0}^\infty2^{-2j\delta}\norm{f}_{L^2}^2\\
&\le C\norm{f}_{L^2}^2,
\end{align*}
which completes the proof.
\end{proof}

We next consider a parametrix construction of $\Op(\chi_\ep)e^{-itH}\Op(\chi_\ep^*)$. Let us first make the following  ansatz: 
$$
v(t,x)=\frac{1}{(2\pi)^d}\int e^{i(\Psi(t,x,\xi)-y\cdot\xi)}b^N(t,x,\xi)f(y)dyd\xi,
$$
where $b^N=\sum_{j=0}^{N-1} b_j$. In order to  approximately solve the Schr\"odinger equation
$$
i\partial_t v(t)=Hv(t);\quad v|_{t=0}=\Op(\chi_\ep^*)\varphi,
$$ 
the phase function $\Psi$ and the amplitude $b^N$ should satisfy the following Hamilton-Jacobi equation and transport equations, respectively:
\begin{align}
\label{Hamilton_Jacobi}
\partial_t\Psi+p(x,\partial_x \Psi)=0;\quad \Psi|_{t=0}=x\cdot\xi,
\end{align}
\begin{equation}
\label{transport_WKB}
\left\{
\begin{aligned}
&\partial_t b_{0}+\X\cdot \partial_xb_{0}+\Y b_{0}=0;\quad b_0|_{t=0}=\chi_\ep,\\
&\partial_t b_{j}+\X\cdot \partial_xb_{j}+\Y b_{j}+iKb_{j-1}=0;\quad b_j|_{t=0}=0,\quad 1\le j \le N-1, 
\end{aligned}
\right.
\end{equation}
where $K$ is the kinetic energy part of $H$, a vector field $\X$ and a function $\Y$ are defined by
\begin{align*}
\X_j(t,x,\xi)
&:=(\partial _{\xi_j} p)(x,\partial_x \Psi(t,x,\xi)),\ j=1,...,d,\\
\Y(t,x,\xi)
&:=[k(x,\partial_x)\Psi+p_1(x,\partial_x\Psi)](t,x,\xi).
\end{align*}
Here $p,p_1$ are given by \eqref{symbols_1}. 
We first construct the phase function $\Psi$.

\begin{proposition}					%proposition
\label{proposition_phase_WKB}
Let us fix $\ep>0$ arbitrarily. Then, for sufficiently small $t_\ep>0$, we can construct a smooth and real-valued function $\Psi \in C^\infty((-t_\ep,t_\ep)\times \R^{2d};\R)$ which solves the Hamilton-Jacobi equation \eqref{Hamilton_Jacobi} for $(x,\xi)\in \Omega(\ep)$ and $|t| \le t_\ep$. Moreover, for all $\alpha,\beta\in \Z^d_+$, $x,\xi \in \R^d$ and $|t|\le t_\ep$,
\begin{align}
\label{proposition_phase_WKB_1}
|\dderiv{x}{\xi}{\alpha}{\beta}(\Psi(t,x,\xi)-x\cdot\xi+tp(x,\xi)|\le C_{\alpha\beta\ep}|t|^2\<x\>^{2-|\alpha+\beta|},
\end{align}
where $C_{\alpha\beta\ep}>0$ is independent of $x,\xi$ and $t$. 
\end{proposition}

\begin{proof}
We consider the case when $t \ge0$, and the proof for $t \le0$ is similar. We first define the action integral $\wtilde{\Psi}(t,x,\xi)$ on $[0,t_\ep) \times \Omega(\ep/2)$ by
$$
\wtilde{\Psi}(t,x,\xi):=x\cdot\xi+\int_0^tL(X(s,Y(t,x,\xi),\xi),\Xi(s,Y(t,x,\xi),\xi))ds, 
$$
where $L(x,\xi)=\xi\cdot\partial_\xi p(x,\xi)-p(x,\xi)$ is the Lagrangian associated to $p(x,\xi)$, and $X,\Xi$ and $Y$ are given by Lemma \ref{WKB_2_lemma_1} (2) with $\ep$ replaced by $\ep/2$. The smoothness of $\wtilde{\Psi}(t,x,\xi)$ follows from corresponding properties of $X(t)$, $\Xi(t)$ and $Y(t)$. It is well known that $\wtilde{\Psi}(t,x,\xi)$ solves the Hamilton-Jacobi equation
$$
\partial_t\wtilde{\Psi}(t,x,\xi)+p(x,\partial_x \wtilde{\Psi}(t,x,\xi))=0;\quad \Psi|_{t=0}=x\cdot\xi,
$$
for $(x,\xi)\in \Omega(\ep/2)$, and satisfies
$$
\partial_x\wtilde{\Psi}(t,x,\xi)=\Xi(t,Y(t,x,\xi),\xi),\quad\partial_\xi \wtilde{\Psi}(t,x,\xi)=Y(t,x,\xi).
$$
Lemma \ref{WKB_2_lemma_1} (2) shows that $p(Y(t,x,\xi),\xi)$ is of quadratic type:
$$
|\dderiv{x}{\xi}{\alpha}{\beta}p(Y(t,x,\xi),\xi)|\le C_{\alpha\beta\ep}\<x\>^{2-|\alpha+\beta|},\quad
(t,x,\xi)\in[0,t_\ep)\times\Omega(\ep/2),
$$ 
which, combined with the energy conservation 
$$
p(x,\partial_x\wtilde{\Psi}(t,x,\xi))=p(Y(t,x,\xi),\xi),
$$ 
imply
$$
|\dderiv{x}{\xi}{\alpha}{\beta}(\wtilde{\Psi}(t,x,\xi)-x\cdot\xi)|
\le C_{\alpha\beta\ep}|t|\<x\>^{2-|\alpha+\beta|},\quad 
(t,x,\xi)\in[0,t_\ep)\times\Omega(\ep/2).
$$ 
We similarly obtain, for $(t,x,\xi)\in[0,t_\ep)\times\Omega(\ep/2)$,
\begin{align*}
&|p(x,\partial_x\wtilde{\Psi}(t,x,\xi))-p(x,\xi)|\\
&=\left|\left(\partial_x\wtilde{\Psi}(t,x,\xi))-\xi\right)\cdot\int_0^1(\partial_\xi p)(x,\theta\partial_x\wtilde{\Psi}(t,x,\xi))+(1-\theta)\xi)d\theta\right|\\
&\le C_\ep |t|\<x\>^2,
\end{align*}
and, more generally, 
\begin{align*}
|\dderiv{x}{\xi}{\alpha}{\beta}(p(x,\partial_x\wtilde{\Psi}(t,x,\xi))-p(x,\xi))|
&\le C_{\alpha\beta\ep}
|t| \<x\>^{2-|\alpha+\beta|}.
\end{align*}
Therefore, integrating the Hamilton-Jacobi equation with respect to $t$, we have
$$
|\dderiv{x}{\xi}{\alpha}{\beta}\left(\wtilde{\Psi}(t,x,\xi)-x\cdot\xi+tp(x,\xi)\right)|
\le C_{\alpha\beta\ep} |t|^2\<x\>^{2-|\alpha+\beta|}.
$$
Finally, choosing a cut-off function $\rho \in S(1,g)$ so that $0 \le \rho\le 1$, $\rho \equiv 1$ on $\Omega(\ep)$ and $\supp \rho \subset \Omega(\ep/2)$, we define $$\Psi(t,x,\xi):=x\cdot\xi-tp(x,\xi)+\rho(x,\xi)(\wtilde{\Psi}(t,x,\xi)-x\cdot\xi+tp(x,\xi)).$$ $\Psi(t,x,\xi)$ clearly satisfies the statement of Proposition \ref{proposition_phase_WKB}. 
\end{proof}

Using the phase function constructed in Proposition \ref{proposition_phase_WKB}, we can define the FIO, $J(\Psi,a):\S\to \S'$ by
$$
J(\Psi,a)f(x)=\frac{1}{(2\pi)^d}\int e^{i(\Psi(t,x,\xi)-y\cdot\xi)}a(x,\xi)f(y)dyd\xi,\quad f \in \S(\R^d),
$$
where $a\in S(1,g)$. Moreover, we have the following.

%{lemma}
\begin{lemma}
\label{WKB_lemma_1}
Let $t_\ep>0$ be small enough. Then, for any bounded family of symbols $\{a(t);|t|<t_\ep\} \subset S(1,g)$, $J(\Psi,a)$ is bounded on $L^2(\R^d)$ uniformly with respect to $|t|<t_\ep$:
$$
\sup_{|t| \le t_\ep}\norm{J(\Psi,a)}_{L^2\to L^2} \le C_\ep.
$$
\end{lemma}

\begin{proof}[Proof] 
For sufficiently small $t_\ep>0$, the estimates \eqref{proposition_phase_WKB_1} imply
$$
|(\partial_\xi\otimes\partial_x \Psi)(t,x,\xi)-\Id|\le C_\ep t_\ep<1/2,\quad |\dderiv{x}{\xi}{\alpha}{\beta}\Psi(t,x,\xi)|\le C_{\alpha\beta\ep}\ \text{for}\ |\alpha+\beta|\ge2,
$$
uniformly with respect to $(t,x,\xi)\in (-t_\ep,t_\ep)\times \R^{2d}$. Therefore, the assertion is a consequence of the standard $L^2$-boundedness of FIO, or equivalently Kuranishi's trick and the $L^2$-boundedness of PDO (see, \emph{e.g.}, \cite{Robert2} or \cite[Lemma 4.2]{Mizutani2}).\end{proof}

We next construct the amplitude.

%proposition
\begin{proposition}					
\label{proposition_amplitude_WKB}
Let $\Psi(t,x,\xi)$ be as in Proposition \ref{proposition_phase_WKB} with $\ep$ replaced by $\ep/3$. Then, for any integer $N\ge0$, there exist families of symbols $\{b_j(t,\cdot,\cdot);|t|<t_\ep\}\subset S(\<x\>^{-j}\<\xi\>^{-j},g)$, $j=0,1,2,...,N-1$, such that $\supp b_j(t,\cdot,\cdot)\subset \Omega(\ep/2)$ and $b_j$ solve the transport equations \eqref{transport_WKB}.
\end{proposition}

\begin{proof}
We consider the case $t\ge0$ only. Symbols $b_{j}$ can be constructed by a standard method of characteristics, along the flow generated by $\X(t,x,\xi)$, as follows. At first note that Assumption \ref{assumption_A} and \eqref{proposition_phase_WKB_1} imply that \begin{align}
\label{proposition_amplitude_WKB_1_1}
|\dderiv{x}{\xi}{\alpha}{\beta}\X(t,x,\xi)| &\le C_{\alpha\beta\ep}\<x\>^{1-|\alpha+\beta|},\\
\label{proposition_amplitude_WKB_1_2}
|\dderiv{x}{\xi}{\alpha}{\beta}\Y(t,x,\xi)| &\le C_{\alpha\beta\ep}\<x\>^{-|\alpha+\beta|},
\end{align}
uniformly with respect to $0\le t\le t_\ep$ and $(x,\xi) \in \Omega(\ep/3)$. 
For all $0 \le s,t \le t_\ep$, we consider the solution to the following ODE:
$$
\partial_t z(t,s,x,\xi)=\X(t,z(t,s,x,\xi),\xi);\quad z(s,s)=x.
$$
We learn by \eqref{proposition_amplitude_WKB_1_1} and a same argument as that in the proof of Lemma \ref{WKB_1_lemma_3} that $z(t,s)$ is well defined for $0 \le s,t \le t_\ep$ and $(x,\xi)\in\Omega(\ep/3)$, and that
\begin{equation}
\begin{aligned}
\label{WKB_theorem_proof_1}
|\dderiv{x}{\xi}{\alpha}{\beta}(z(t,s,x,\xi)-x)| 
\le C_{\alpha\beta\ep}t_\ep \<x\>^{1-|\alpha+\beta|},\quad (x,\xi)\in\Omega(\ep/3). 
\end{aligned}
\end{equation}
Then, $b_{j}(t)$ are defined inductively by
\begin{equation}
\begin{aligned}
\nonumber
b_{0}(t,x,\xi)&=\chi_\ep^*(z(0,t,x,\xi),\xi)\exp\left(\int_0^t\Y(s,z(s,t,x,\xi),\xi)ds\right),\\
b_{j}(t,x,\xi)&=-\int_0^t (iKb_{j-1})(s,z(s,t,x,\xi),\xi) \exp\left(\int_u^t\Y(u,z(u,t,x,\xi),\xi)du\right)ds. 
\end{aligned}
\end{equation}
Since $\supp \chi_\ep^* \subset \Omega(\ep)$, by \eqref{WKB_theorem_proof_1} and a same argument as that in the proof of Lemma \ref{WKB_2_lemma_1} (1), we see that $b_{j}(t,x,\xi)$ is smooth with respect to $(x,\xi)$ and that $\dderiv{x}{\xi}{\alpha}{\beta}b_{j}(t,x,\xi)$ are supported in $\Omega(\ep/2)$ for all $0\le t\le t_\ep$. Thus, if we extend $b_{j}$ on $\R^{2d}$ so that $b_{j}(t,x,\xi)=0$ outside $\Omega(\ep/2)$, then $b_{j}$ is still smooth in $(x,\xi)$. Furthermore, we learn by \eqref{proposition_amplitude_WKB_1_2} and \eqref{WKB_theorem_proof_1} that $\{b_{j}(t,\cdot,\cdot);t \in [0,t_\ep],\ 0 \le j \le N-1\}$ is a bounded set in $S(\<x\>^{-j}\<\xi\>^{-j},g)$. Finally, a standard Hamilton-Jacobi theory shows that $b_j(t)$ solve the transport equations \eqref{transport_WKB}.
\end{proof}

We now state the main result in this section.

%theorem
\begin{theorem}						
\label{WKB_theorem_1}
Let us fix $\ep>0$ arbitrarily. Then, for sufficiently small $t_\ep>0$, any nonnegative integer $N\ge0$ and any symbol $\chi_\ep \in S(1,g)$ supported in $\Omega(\ep)$, 
we can find a bounded family of symbols 
$
\{a^N(t,\cdot,\cdot);|t|<t_\ep\}\subset S(1,g)
$
such that $\Op(\chi_\ep)e^{-itH}\Op(\chi_\ep)^*$ can be brought to the form
\begin{align*}
\Op(\chi_\ep)e^{-itH}\Op(\chi_\ep)^*=J(\Psi,a^N)+Q(t,N),
\end{align*}
where $J(\Psi,a^N)$ is the FIO with the phase $\Psi(t,x,\xi)$ constructed in Proposition \ref{proposition_phase_WKB} with $\ep$ replaced by $\ep/3$. The distribution kernel of $J(\Psi,a^N)$, which we denote by $K_{\Psi,a^N}(t,x,y)$, satisfies the dispersive estimate:
$$
|K_{\Psi,a^N}(t,x,y)|\le C_{N,\ep}|t|^{-d/2},\quad 0<|t|< t_\ep,\ x,\xi\in \R^d.
$$
Moreover, for any $\gamma\ge0$ with $N>2\gamma$, the remainder $Q(t,N)$ satisfies 
\begin{align}
\label{WKB_theorem_1_1}
\norm{\<D\>^\gamma Q(t,N)\<D\>^{\gamma}}_{L^2 \to L^2}
\le C_{N\gamma\ep} |t|,\quad |t| <t_\ep.
\end{align}
In particular, if we choose $N\ge d+1$, then the distribution kernel of $Q(t,N)$ is uniformly bounded in $\R^{2d}$ with respect to $|t|<t_\ep$, and hence
$$
\norm{\Op(\chi_\ep)e^{-itH}\Op(\chi_\ep)^*}_{L^1 \to L^\infty} \le C_\ep|t|^{-d/2},
\quad 0<|t|<t_\ep.
$$
\end{theorem}

\begin{proof}
We consider the case when $t \ge0$ and the proof for the opposite case is similar. By virtue of Lemma \ref{lemma_WKB}, we may replace $\Op(\chi_\ep)^*$ by $\Op(\chi_\ep^*)$ for some $ \chi_\ep^*\in S(1,g)$  supported in $\Omega(\ep)$, without loss of generality.
Let $b^N=\sum_{j=0}^{N-1}b_j$ with $b_j$ constructed in Proposition \ref{proposition_amplitude_WKB}. Since $J(\Psi,b^N)|_{t=0}=\Op(\chi_\ep^*)$, we have the Duhamel formula
\begin{align*}
\Op(\chi_\ep)e^{-itH}\Op(\chi_\ep^*)
=\Op(\chi_\ep)J(\Psi,b^N)-i \int_0^t \Op(\chi_\ep) e^{-i(t-s)H}(D_t+H)J(\Psi,b^N)|_{t=s}ds.
\end{align*}
\textbf{Estimates on the remainder.} It suffices to show that
$$
\sup_{|t|<t_\ep}\norm{\<D\>^\gamma \Op(\chi_\ep) e^{-itH}(D_t+H)J(\Psi,b^N)\<D\>^\gamma}_{L^2\to L^2} \le C_{N\gamma\ep}.
$$
Since $\Psi,b_j$ solve the Hamilton-Jacobi equation \eqref{Hamilton_Jacobi} and transport equations \eqref{transport_WKB}, respectively, a direct computation yields 
$$
e^{-i\Psi(t,x,\xi)}(D_t+H)\left(e^{i\Psi(t,x,\xi)}\sum_{j=0}^{N-1} b_j(t,x,\xi)\right)
=r_{N}(t,x,\xi),
$$
with some $\{r_{N}(t,\cdot,\cdot);0\le t\le t_\ep\} \subset S(\<x\>^{-N}\<\xi\>^{-N},g)$. In particular, 
$$
(D_t+H)J(\Psi,b^N)=J(\Psi,r_N).
$$
A standard $L^2$-boundedness of FIO then implies 
$$
\sup_{|t|< t_\ep}\norm{\<x\>^{\gamma+\delta}\<D\>^{\gamma+\delta}J(\Psi,r_N)\<D\>^\gamma}_{L^2 \to L^2} \le C_{N\gamma\delta}, 
$$
for any $\gamma,\delta\ge0$ with $2\gamma+\delta\le N$. Since we already proved in the proof of Lemma \ref{lemma_WKB} that 
$$
\sup_{|t|\le t_\ep}\norm{\<D\>^\gamma\Op(\chi_\ep)e^{-itH}\<D\>^{-\gamma-\delta}\<x\>^{-\gamma-\delta}}_{L^2 \to L^2} \le C_{\gamma\delta}, 
$$
we obtain the desired estimate. 
\\\textbf{Dispersive estimates.} By the composition formula of PDO and FIO (cf. \cite{Robert2}), $\Op(\chi_\ep)J(\Psi,b^N)$ is also a FIO with the same phase $\Psi$ and the amplitude
$$
a^N(t,x,\xi)=\frac{1}{(2\pi)^d}\int e^{iy\cdot\eta}\chi_\ep(x,\eta+\wtilde{\Xi}(t,x,y,\xi))b^N(t,x+y,\xi)dyd\eta,
$$
where $\wtilde{\Xi}(t,x,y,\xi)=\int_0^1 (\partial_x \Psi)(t,y+\lambda(x-y),\xi)d\lambda$. By virtue of \eqref{proposition_phase_WKB_1}, $\wtilde{\Xi}$ satisfies
$$
|\dderiv{x}{y}{\alpha}{\alpha'}\partial_\xi^\beta(\wtilde{\Xi}(t,x,y,\xi)-\xi)| \le C_{\alpha\alpha'\beta}|t|,\quad |\alpha+\alpha'+\beta|\ge 1. 
$$
Combining with the fact that $\chi_\ep,b^N\in S(1,g)$, $\supp \chi_\ep\subset \Omega(\ep)$ and $\supp b^N(t,\cdot,\cdot) \subset \Omega(\ep/2)$, we see that $\{a^N; 0\le t < t_\ep\}$ is bounded in $S(1,g)$. The distribution kernel of $J(\Psi,a^N)$ is given by
$$K_{\Psi,a^N}(t,x,y)=\frac{1}{(2\pi)^d}\int e^{i(\Psi(t,x,\xi)-y\cdot\xi)}a^N(t,x,\xi)d\xi. $$
By virtue of Proposition \ref{proposition_phase_WKB}, we have
\begin{align*}
&\sup_{|t|\le t_\ep}|\dderiv{x}{y}{\alpha}{\beta}\partial^\gamma_\xi(\Psi(t,x,\xi)-y\cdot\xi)|\le C_{\alpha\beta\gamma},\quad |\alpha+\beta+\gamma|\ge 2,\\
&\partial_\xi^2\Psi(t,x,\xi)=-t(g^{jk}(x))_{j,k}+O(t^2),\quad |t|\to 0.
\end{align*}
As a consequence, since $g^{jk}(x)$ is uniformly elliptic, the phase function $\Psi(t,x,\xi)-y\cdot\xi$ has a unique non-degenerate critical point for all $|t|<t_\ep$ and we can apply the stationary phase method to $K_{\Psi,a^N}(t,x,y)$, provided that $t_\ep>0$ is small enough. Therefore, 
$$|K_{\Psi,a^N}(t,x,y)|\le C|t|^{-d/2},\quad 0<|t|\le t_\ep,\ x,\xi\in \R^d, $$
which completes the proof.
\end{proof}

%%%%%%%%%%%%%%%%%%%%%%%%%%%%%%%%%%%%%%%%%%%%%%%%%%%		Proof of Theorem \ref{theorem_2}		      %%%%%%%%%%%%%%%%%%%%%%%%%%%%%%%%%%%%%%%%%%%%%%%%%%%

\section{Proof of Theorem \ref{theorem_2}}
\label{section_proof_theorem_2}
Suppose that $H$ satisfies Assumption \ref{assumption_A} with $\mu\ge0$. In this section we give the proof of Theorem \ref{theorem_2}. In view of Corollary \ref{corollary_Littlewood_Paley}, \eqref{theorem_2_1} is a consequence of the following proposition.

%proposition
\begin{proposition}
\label{proposition_proof_theorem_2_1}
For any symbol $a\in C_0^\infty(\R^{2d})$ and $T>0$, 
$$
\norm{\Op_h(a)e^{-itH}\varphi}_{L^p([-T,T];L^q(\R^d))} \le C_Th^{-1/p}\norm{\varphi}_{L^2(\R^d)},
$$
uniformly with respect to $h\in(0,1]$, provided that $(p,q)$ satisfies \eqref{admissible}. 
\end{proposition}

\begin{proof}
This proposition follows from the standard WKB approximation for $e^{-itH}\Op_h(a)$ up to time scales of order $1/h$. The proof is essentially same as that in the case for the Laplace-Beltrami operator on compact manifolds without boundaries (see, \cite[Section 2]{BGT}), and we omit details.
\end{proof}

Using this proposition, we have the semiclassical Strichartz estimates with inhomogeneous error terms:

%proposition
\begin{proposition}					
\label{proposition_proof_theorem_2_2}
Let $a\in C_0^\infty(\R^{2d})$. Then, for any $T>0$ and any $(p,q)$ satisfying the admissible condition \eqref{admissible},
\begin{align*}
\norm{\Op_h(a)e^{-itH}\varphi}_{L^p([-T,T];L^q(\R^d))}
&\le C_T\norm{\Op_h(a)\varphi}_{L^2(\R^d)}+C_Th\norm{\varphi}_{L^2(\R^d)}\\
&+Ch^{-1/2}\norm{\Op_h(a)e^{-itH}\varphi}_{L^2([-T,T];L^2(\R^d))}\\
&+Ch^{1/2}\norm{[\Op_h(a),H]e^{-itH}\varphi}_{L^2([-T,T];L^2(\R^d))},
\end{align*}
uniformly with respect to $h\in(0,1]$. 
\end{proposition}

This proposition has been proved by \cite{Bouclet_Tzvetkov_1} for the case with $V,A\equiv0$. We give a refinement of this proposition with its proof in Section \ref{section_proof_theorem_3}.

Next, we shall prove that if $k(x,\xi)$ satisfies the nontrapping condition \eqref{nontrapping}, then the missing $1/p$ derivative can be recovered. We first recall the local smoothing effects for Schr\"odinger operators proved by Doi \cite{Doi}. For  any $s \in \R$, we set 
$
\B^s:=\{ f \in L^2(\R^d) ; \<x\>^s f, \<D\>^s f \in L^2(\R^d)\}
$. 
Define a symbol $e_s(x,\xi)$ by
$$
e_s(x,\xi):=(k(x,\xi)+|x|^2+L(s))^{s/2}\in S((1+|x|+|\xi|)^s,g), 
$$
where $L(s)>1$ is a large constant depending on $s$. We denote by $E_s$ its Weyl quantization:
$$
E_sf(x)=\Op^w(e_s)f(x)=\frac{1}{(2\pi)^{d}}\int e^{i(x-y)\cdot\xi}e_s\left(\frac{x+y}{2},\xi\right)f(y)dyd\xi.
$$
Then, for any $s \in \R$, there exists $L(s)>0$ such that $E_s$ is a homeomorphism from $\B^{r+s}$ to $\B^r$ for all $r \in \R$, and $(E_s)^{-1}$ is still a Weyl quantization of a symbol in $S((1+|x|+|\xi|)^{-s},g)$ (see, \cite[Lemma 4.1]{Doi}). 

%proposition
\begin{proposition}[The local smoothing effects \cite{Doi}]
\label{local_smoothing_effect_1}
Suppose that $k(x,\xi)$ satisfies the nontrapping condition \eqref{nontrapping} and that Assumption \ref{assumption_B}. Then, for any $T>0$ and $\sigma>0$, there exists $C_{T,\sigma}>0$ such that
\begin{align}
\label{local_smoothing_effect_1_1}
\norm{\<x\>^{-1/2-\sigma}E_{1/2}e^{-itH}\varphi}_{{L^2([-T,T];L^2(\R^d))}} \le C_{T,\sigma} \norm{\varphi}_{L^2(\R^d)}.
\end{align}
\end{proposition}

%remark
\begin{remark}\eqref{local_smoothing_effect_1_1} implies a standard local smoothing effect:
\begin{align}
\label{local_smoothing_effect_2}
\norm{\<x\>^{-1/2-\sigma}\<D\>^{1/2}e^{-itH}\varphi}_{L^2([-T,T];L^2(\R^d))} \le C_{T,\sigma}\norm{\varphi}_{L^2(\R^d)}.
\end{align}
 Indeed, we compute
\begin{align*}
&\<x\>^{-1/2-\sigma}\<D\>^{1/2}\\
&=\<D\>^{1/2}\<x\>^{-1/2-\sigma}+[\<D\>^{1/2},\<x\>^{-1/2-\sigma}]\\
&=\<D\>^{1/2}(E_{1/2})^{-1}E_{1/2}\<x\>^{-1/2-\sigma}+[\<D\>^{1/2},\<x\>^{-1/2-\sigma}]\\
&=\<D\>^{1/2}(E_{1/2})^{-1}\left(\<x\>^{-1/2-\sigma}E_{1/2}+[E_{1/2},\<x\>^{-1/2-\sigma}]\right)+[\<D\>^{1/2},\<x\>^{-1/2-\sigma}].
\end{align*}
It is easy to see that $\<D\>^{1/2}(E_{1/2})^{-1}$, $[E_{1/2},\<x\>^{-1/2-\sigma}]$ and $[\<D\>^{1/2},\<x\>^{-1/2-\sigma}]$ are bounded on $L^2(\R^d)$ since their symbols belong to $S(1,g)$. Therefore, \eqref{local_smoothing_effect_1_1} implies \eqref{local_smoothing_effect_2}. 
\end{remark}

\begin{proof}[Proof of \eqref{theorem_2_2} of Theorem \ref{theorem_2}]
It is clear that \eqref{theorem_2_2} follows from Proposition \ref{proposition_proof_theorem_2_2}, \eqref{local_smoothing_effect_2} and Corollary \ref{corollary_Littlewood_Paley}, since $a$ is compactly supported with respect to $x$ and $\{a,p\}\in S(\<\xi\>,g)$, where $p=p(x,\xi)$. 
\end{proof}

%%%%%%%%%%%%%%%%%%%%%%%%%%%%%%%%%%%%%%%%%%%%%%%%%			Strichartz estimates with loss without asymptotic flatness    %%%%%%%%%%%%%%%%%%%%%%%%%%%%%%%%%%%%%%%%%%%%

\section{Near sharp Strichartz estimates without asymptotic flatness}
\label{section_proof_theorem_3}
This section is devoted to prove Theorem \ref{theorem_3}. We may assume $\mu=0$ without loss of generality. We begin with the following proposition. 
%proposition
\begin{proposition}
\label{proof_theorem_3_proposition_1_1}
Let $I \Subset(0,\infty)$ be a relatively compact open interval and $C_0>1$. Then, there exist $\delta_0,h_0>0$ such that for any $0<\delta\le\delta_0$, $0<h\le h_0$, $1\le R\le 1/h$ and any symbol $a_h\in S(1,g)$ supported in $\{(x,\xi);R<|x|<C_0/h,\ |\xi|\in I\}$, we have
\begin{align}
\label{proof_theorem_3_proposition_1}
\norm{\Op_h(a_h)e^{-itH}\Op_h(a_h)^*}_{L^1 \to L^\infty} \le C_\delta|t|^{-d/2},\quad 0<|t|<\delta hR,
\end{align}
where $C_\delta>0$ may be taken uniformly with respect to $h$ and $R$.
\end{proposition}

\begin{remark}
When $|t|>0$ in \eqref{proof_theorem_3_proposition_1} is small and independent of $R$, Proposition \ref{proof_theorem_3_proposition_1} is well-known and the proof is given by the standard method of the short-time WKB approximation for $e^{-itH^h/h}\Op_h(a_h)^*$ (see, \emph{e.g.}, \cite{BGT}). 
\end{remark}

For $h\in(0,1]$, $R\ge1$, an open interval $I\Subset(0,\infty)$ and $C_0>1$, we set 
$$
\Gamma(R,h,I):=\{(x,\xi)\in\R^{2d};R<|x|<C_0/h,\ |\xi|\in I\}.
$$ 
Proposition \ref{proof_theorem_3_proposition_1} is a consequence of the same argument as in the proof of Proposition \ref{theorem_dispersive_IK} and the following proposition:

%theorem
\begin{proposition}
Let $I \Subset I_1 \Subset (0,\infty)$ and $C_0>1$. Then, there exist $\delta_0,h_0>0$ such that the followings hold for any $0<\delta\le\delta_0$, $0<h\le h_0$ and $1\le R \le C_0/h$.\\
\emph{(1)} There exists $\Phi_h(t,x,\xi)\in C^\infty((-\delta R,\delta R)\times\R^{2d})$ such that $\Phi_h$ solves  the following Hamilton-Jacobi equation:
\begin{equation}\left\{\begin{aligned}
\label{proof_theorem_3_proposition_1_1_10}
\partial_t \Phi_h(t,x,\xi)
&=-p_h(x,\partial_x\Phi_h(t,x,\xi)),\quad |t|<\delta R,\ (x,\xi) \in \Gamma(R/2,h/2,I_1), \\
\Phi_h(0,x,\xi)
&=x\cdot\xi,\quad 
(x,\xi)\in \Gamma(R/2,h/2,I_1).
\end{aligned}\right.\end{equation}
Furthermore, we have
\begin{align}
\label{proof_theorem_3_proposition_1_2}
|\dderiv{x}{\xi}{\alpha}{\beta}\left(\Phi_h(t,x,\xi)-x\cdot\xi+tp_h(x,\xi)\right)|
\le C_{\alpha\beta} R^{-|\alpha|}h|t|^2,\quad \alpha,\beta\in \Z^d_+,
\end{align} 
uniformly with respect to $x,\xi\in\R^d$, $h\in (0,h_0]$, $0 \le R \le C_0/h$ and $|t|<\delta R$. \\
\emph{(2)} For any $a_h \in S(1,g)$ with $\supp a_h \subset \Gamma(R,h,I)$ and any integer $N\ge0$, we can find $b_{h}^N(t,\cdot,\cdot)\in S(1,g)$ such that 
\begin{align*}
e^{-it\wtilde{H}^h/h}\Op_h(a_h)^*=J_h(\Phi_h,b_h^N)+Q_{{\WKB}}(t,h,N),
\end{align*}
where $J_h(\Phi_h,b_h^N)$ is the $h$-FIO with the phase function $\Phi_h$ and the amplitude $b_h^N$, and its distribution kernel satisfies
\begin{align}
\label{dispersive_2}
|K_{{\WKB}}(t,h,x,y)| \le C |th|^{-d/2},\quad h \in (0,h_0],\ 0<|t| \le \delta R,\ x,\xi \in \R^d.
\end{align}
Moreover the remainder $Q_{{\WKB}}(t,h,N)$ satisfies
\begin{align*}
\norm{\<D\>^sQ_{{\WKB}}(t,h,N)\<D\>^s}_{L^2 \to L^2}
\le C_{N,s} h^{N-2s}|t|,\quad h \in (0,h_0],\ |t| \le \delta R.
\end{align*}
\end{proposition}

\begin{proof}[Sketch of the proof]
The proof is similar to that of Theorem \ref{WKB_theorem_1} and, in particular, the proof of the second claim is completely same. Thus, we give only the outline of the construction of $\Phi_h$. We may assume $C_0=1$ without loss of generality. Let us denote by $(X_h,\Xi_h)$ the Hamilton flow generated by $p_h$. To construct the phase function, the most important step is to study the inverse map of $(x,\xi)\mapsto (X_h(t,x,\xi),\xi)$. Choose an open interval $\wtilde{I}_1$ so that $I_1\Subset \wtilde{I}_1\Subset(0,\infty)$. The following bounds have been proved by \cite{Mizutani2}:
$$
|\dderiv{x}{\xi}{\alpha}{\beta}(X_h(t,x,\xi)-x)|+\<x\>|\dderiv{x}{\xi}{\alpha}{\beta}(\Xi_h(t,x,\xi)-\xi)|\le C_{\alpha\beta}\<x\>^{-|\alpha|}|t|,
$$
for $(x,\xi)\in\Gamma(R/3,h/3,\wtilde{I}_1)$ and $|t|\le \delta R$. For sufficiently small $\delta>0$ and for any fixed $|t|\le \delta R$, the above estimates imply
$$
|\partial_xX_h(t)-\Id| 
\le CR^{-1}|t|\le C\delta <1/2.
$$
By the same argument as that in the proof of Lemma \ref{WKB_2_lemma_1}, the map $(x,\xi) \mapsto (X_h(t,x,\xi),\xi)$ is a diffeomorphism from $\Gamma(R/3,h/3,\wtilde{I}_1)$ onto its range and the corresponding inverse $(x,\xi)\mapsto (Y_h(t,x,\xi),\xi)$ is well-defined for $|t|<\delta R$ and $(x,\xi)\in \Gamma(R/2,h/2,I_1)$. Moreover, $Y_h(t)$ satisfies the same estimates as that for $X_h(t)$:
$$
|\dderiv{x}{\xi}{\alpha}{\beta}(Y_h(t,x,\xi)-x)|\le C_{\alpha\beta}\<x\>^{-|\alpha|}|t|,\quad |t|<\delta R,\ (x,\xi)\in \Gamma(R/2,h/2,I_1). 
$$
We now define $\Phi_h$ by 
$$
\Phi_h(t,x,\xi):=x\cdot\xi+\int_0^tL_h(X_h(s,Y(t,x,\xi),\xi),\Xi(s,y(t,x,\xi),\xi))ds, 
$$
where $L_h=\xi\cdot\partial_\xi p_h-p_h$. By the standard Hamilton-Jacobi theory, $\Phi_h$ solves \eqref{proof_theorem_3_proposition_1_1_10}. Moreover, using the energy conservation $p_h(x,\partial_x\Phi_h(t))=p_h(Y_h(t),\xi)$ and the above estimates on $X_h,\Xi_h$ and $Y_h$, we see that
\begin{align*}
|p_h(x,\partial_x\Phi_h(t))-p_h(x,\xi)|
&=|p_h(Y_h(t),\xi)-p_h(x,\xi)|
\\
&\le |Y_h(t)-x|\left|\int_0^\lambda (\partial_x p_h)(\lambda Y_h(t)-(1-\lambda)x,\xi)d\lambda\right|\\
&\le C|y(t)-x|(h+h^2\<x\>^2)\\
&\le Ch|t|,
\end{align*}
and that
$$
|\dderiv{x}{\xi}{\alpha}{\beta}(p_h(x,\partial_x\Phi_h)-p_h(x,\xi))| \le C_{\alpha\beta}\<x\>^{-|\alpha|}h|t|. 
$$
Using these estimates, we can check that $\Phi_h$ satisfies \eqref{proof_theorem_3_proposition_1_2}. Finally, we extend $\Phi_h$ to the whole space so that $\Phi_h(t,x,\xi)=x\cdot\xi-tp_h(x,\xi)$ outside $\Gamma(R/3,h/3,\wtilde{I}_1)$.
\end{proof}

Using Proposition \ref{proof_theorem_3_proposition_1_1}, we obtain a refinement of Proposition \ref{proposition_proof_theorem_2_2}:

%proposition
\begin{proposition}
\label{proof_theorem_3_proposition_2}
Let $0<R\le 1/h$ and let $a_h\in S(1,g)$ be supported in $\{(x,\xi);R<|x|<1/h,\ |\xi|\in I\}$. Then, for any $T>0$ and $(p,q)$ satisfying the admissible condition \eqref{admissible},
\begin{align*}
\norm{\Op_h(a_h)e^{-itH}\varphi}_{L^p([-T,T];L^q(\R^d))}
&\le C_T\norm{\Op_h(a_h)\varphi}_{L^2(\R^d)}+C_Th\norm{\varphi}_{L^2(\R^d)}\\
&+C_T(hR)^{-1/2}\norm{\Op_h(a_h)e^{-itH}\varphi}_{L^2([-T,T];L^2(\R^d))}\\
&+C_T(hR)^{1/2}\norm{[H,\Op_h(a_h)]e^{-itH}\varphi}_{L^2([-T,T];L^2(\R^d))},
\end{align*}
uniformly with respect to $h\in(0,h_0]$. 
\end{proposition}

\begin{proof}
The proof is similar to that of \cite[Proposition 5.4]{Bouclet_Tzvetkov_1}. By time reversal invariance we can restrict our considerations to the interval $[0,T]$. We may assume $T\ge hR$ without loss of generality and split $[0,T]$ as follows: $[0,T]=J_0\cup J_1\cup\dots\cup J_N$, where $J_j=[jhR,(j+1)hR]$, $0\le j\le N-1$, and $J_N=[T-\delta hR,T]$. For $j=0$, we have the Duhamel formula
$$
\Op_h(a_h)e^{-itH}=e^{-itH}\Op_h(a_h)-i\int_0^te^{-i(t-s)H}[\Op_h(a_h),H]e^{-isH}ds,\quad t\in J_0.
$$
We here choose $b_h\in S(1,g)$ so that $b_h\equiv 1$ on $\supp a$ and $b_h$ is supported in a sufficiently small neighborhood of $\supp a_h$. By Proposition \ref{proof_theorem_3_proposition_1_1}, $\Op_h(b_h)e^{-i(t-s)H}\Op_h(b_h)^*$ satisfies dispersive estimates \eqref{proof_theorem_3_proposition_1} for $0<|t-s|<\delta hR$ with some $\delta>0$ small enough. Using the Keel-Tao theorem \cite{Keel_Tao} and the unitarity of $e^{-itH}$, we then learn that for any interval $J_R$ of size $|J_R|\le 2hR$, the following homogeneous and inhomogeneous Strichartz estimates hold uniformly with respect to $h\in(0,h_0]$:
\begin{align}
\label{5001}
\norm{\Op_h(b_h)e^{-itH}\varphi}_{L^p(J_R;L^q(\R^d))} 
\le C\norm{\varphi}_{L^2(\R^d)},
\end{align}
\begin{equation}
\begin{aligned}
\label{5002}
\bignorm{\int_0^tF(s\in J_R)\Op_h(b_h)e^{-i(t-s)H}\Op_h(b_h)^*g(s)ds}_{L^p(J_R;L^q(\R^d))}
\le C\norm{g}_{L^1(J_R;L^2(\R^d))},
\end{aligned}
\end{equation}
where $F(s\in J_R)$ is the characteristic function of $J_R$ and $(p,q)$ satisfies the admissible condition \eqref{admissible}.
 On the other hand, using the expansions \eqref{asymptotic_expansion_composition} and \eqref{asymptotic_expansion_adjoint}, we see  that for any $M\ge0$, 
\begin{align*}
\Op_h(a_h)
&=\Op_h(b_h)\Op_h(a_h)+h^M\Op_h(r_{1,h})\\
&=\Op_h(b_h)^*\Op_h(a_h)+h^M\Op_h(r_{2,h}),\\
[\Op_h(a_h),H]&=\Op_h(b_h)^*[\Op_h(a_h),H]+h^M\Op_h(r_{3,h}),
\end{align*}
with some $\{r_{l,h}\}_{h\in(0,1]}$, $l=1,2,3$, which are bounded in $S(\<x\>^{-M}\<\xi\>^{-M},g)$. Therefore, we can write
\begin{align*}
\Op_h(a_h)e^{-itH}
&=\Op_h(b_h)e^{-itH}\Op_h(a_h)\\
&-i\int_0^t\Op_h(b_h)e^{-i(t-s)H}\Op_h(b_h)^*[\Op_h(a_h),H]e^{-isH}ds\\
&+Q(t,h,M),
\end{align*}
where the remainder $Q(t,h,M)$ satisfies
$$
\norm{Q(t,h,M)}_{L^2 \to L^q} \le C_Mh^{M-1-d(1/2-1/q)},\quad 2 \le q \le \infty,
$$
uniformly in $h\in(0,1]$. Combining this estimate with \eqref{5001} and \eqref{5002}, we obtain
\begin{align*}
&\norm{\Op_h(a_h)e^{-itH}}_{L^p(J_0;L^q)}\\
&\le C\norm{\Op_h(a_h)\varphi}_{L^2}+Ch\norm{\varphi}_{L^2}
+ C\norm{[\Op_h(a_h),H]e^{-itH}\varphi}_{L^1(J_0;L^2)}\\
&\le C\norm{\Op_h(a_h)\varphi}_{L^2}+Ch\norm{\varphi}_{L^2}
+ C(hR)^{1/2}\norm{[\Op_h(a_h),H]e^{-itH}\varphi}_{L^2(J_0;L^2)}.
\end{align*}
We similarly obtain the same bound for $j=N$:
\begin{align*}
&\norm{\Op_h(a_h)e^{-itH}}_{L^p(J_N;L^q)}\\
&\le C\norm{\Op_h(a_h)\varphi}_{L^2}+Ch\norm{\varphi}_{L^2}
+ C(hR)^{1/2}\norm{[\Op_h(a_h),H]e^{-itH}\varphi}_{L^2(J_N;L^2)}.
\end{align*}
For $j=1,2,...,N-1$, taking $\theta \in C_0^\infty(\R)$ so that $\theta\equiv1$ on $[-1/2,1/2]$ and $\supp \theta \subset[-1,1]$, we set $\theta_j(t)=\theta(t/(hR)-j-1/2))$. It is easy to see that $\theta_j\equiv1$ on $J_j$ and $\supp \theta_j\subset \wtilde{J}_j=J_j+[-hR/2,hR/2]$. We consider $v_j=\theta_j(t)\Op_h(a_h)e^{-itH}\varphi$, which solves
$$
i\partial_t v_j=Hv_j+\theta_j' \Op_h(a_h)e^{-itH}\varphi+\theta_j[\Op_h(a_h),H]e^{-itH}\varphi;\quad v_j|_{t=0}=0.
$$
A same argument as above and the Duhamel formula then imply that, for any $t\in \wtilde{J}_j$ and $M\ge0$, $v_j$ satisfies
\begin{align*}
v_j&=-i\int_0^t \Op_h(b_h)e^{-i(t-s)H}\Op_h(b_h)^*\left(\theta_j'(s) \Op_h(a_h)+\theta_j(s)[\Op_h(a_h),H]\right)e^{-isH}\varphi ds\\
&+\wtilde{Q}(t,h,M),
\end{align*}
where the remainder $\wtilde{Q}(t,h,M)$ satisfies
$$
\norm{\wtilde{Q}(t,h,M)}_{L^2 \to L^q} \le C_Mh^{M-1-d(1/2-1/q)},\quad 2 \le q \le \infty,
$$
uniformly in $h\in(0,1]$ and $t\in \wtilde{J}_j$. Taking $M\ge0$ large enough,  we learn
\begin{align*}
&\norm{v_j}_{L^p(J_j;L^q)}\\
&\le Ch^2\norm{\varphi}_{L^2}
+ C(hR)^{-1}\norm{\Op_h(a_h)e^{-itH}\varphi}_{L^1(\wtilde{J}_j;L^2)}
+ C\norm{[\Op_h(a_h),H]e^{-itH}\varphi}_{L^1(\wtilde{J}_j;L^2)}\\
&\le Ch^2\norm{\varphi}_{L^2}
+ C(hR)^{-1/2}\norm{\Op_h(a_h)e^{-itH}\varphi}_{L^2(\wtilde{J}_j;L^2)}\\
&\ \ \ \ \ \ \ \ \quad\quad\quad\ \!+ C(hR)^{1/2}\norm{[\Op_h(a_h),H]e^{-itH}\varphi}_{L^2(\wtilde{J}_j;L^2)}.
\end{align*}
Summing over $j=0,1,...,N$, since $N \le T/h$ and $p\ge2$, we have the assertion by Minkowski's inequality.
\end{proof}

\begin{proof}[Proof of Theorem \ref{theorem_3}]
In view of Corollary \ref{corollary_Littlewood_Paley}, Theorem \ref{theorem_2} and Proposition \ref{theorem_dispersive_WKB}, it suffices to show that, for any $a_h\in S(1,g)$ with $\supp a_h\in \{(x,\xi);2\le|x|\le 1/h,\ |\xi|\in I\}$ and any $\ep>0$, 
$$
\sum_{h}\norm{\Op_h(a_h)e^{-itH}f(h^2H)\varphi}_{L^p([-T,T];L^q)}^2 \le C_{T,\ep}\norm{\<H\>^{\ep}\varphi}_{L^2}^2.
$$
Let us consider a dyadic partition of unity:
$$
\sum_{1 \le j \le j_h}\chi(2^{-j}x)=1,\quad 2\le |x| \le 1/h,
$$
where $\chi \in C_0^\infty(\R^d)$ with $\supp \chi\subset\{1/2<|x|<2\}$ and $j_h\le [\log(1/h)]+1$. We set $\chi_j(x)=\chi(2^{-j}x)$. Proposition \ref{proof_theorem_3_proposition_2} then implies
\begin{align*}
\norm{\chi_j\Op_h(a_h)e^{-itH}\varphi}_{L^p([-T,T];L^q)}
&\le C_T\norm{\chi_j\Op_h(a_h)\varphi}_{L^2}+C_Th\norm{\varphi}_{L^2}\\
&+C_T(h2^j)^{-1/2}\norm{\chi_j\Op_h(a_h)e^{-itH}\varphi}_{L^2([-T,T];L^2)}\\
&+C_T(h2^j)^{1/2}\norm{[\chi_j\Op_h(a_h),H]e^{-itH}\varphi}_{L^2([-T,T];L^2)}.
\end{align*}
Since $2^{j-1}\le|x|\le 2^{j+1}$ and $|x|\le 1/h$ on $\supp \chi_ja_h$ we have, for any $\ep\ge0$,
\begin{align*}
&(h2^j)^{-1/2}\norm{\chi_j\Op_h(a_h)e^{-itH}\varphi}_{L^2([-T,T];L^2)}\\
&\le C\norm{\chi_j\<x\>^{-1/2-\ep}h^{-1/2-\ep}\Op_h(a_h)e^{-itH}\varphi}_{L^2([-T,T];L^2)}.
\end{align*}
Since $\{\chi_ja_h,p\}\in S(\<x\>^{-1}\<\xi\>,g)$, we similarly obtain
\begin{align*}
&(h2^j)^{1/2}\norm{\chi_j[\Op_h(a_h),H]e^{-itH}\varphi}_{L^2([-T,T];L^2)}\\
&\le \norm{\wtilde{\chi}_j\<x\>^{-1/2-\ep}h^{-1/2-\ep}\Op_h(b_h)e^{-itH}\varphi}_{L^2([-T,T];L^2)}+C_Th\norm{\varphi}_{L^2},
\end{align*}
where $\wtilde{\chi}_j(x)=\wtilde{\chi}(2^{-j}x)$ for some $\wtilde{\chi}\in C_0^\infty(\R^d)$ satisfying $\wtilde{\chi}\equiv1$ on $[1/2,2]$ and $\supp \wtilde{\chi}\subset[1/4,4]$, and $b_h\in S(1,g)$ is supported in a neighborhood of $\supp a_h$ so that $b_h\equiv1$ on $\supp a_h$. Summing over $1\le j\le j_h$ and using the local smoothing effect \eqref{local_smoothing_effect_2}, since $p,q\ge2$, we obtain
\begin{align*}
&\norm{\Op_h(a_h)e^{-itH}\varphi}_{L^p([-T,T];L^q)}^2\\
&\le \sum_{1\le j\le j_h}\norm{\chi_j\Op_h(a_h)e^{-itH}\varphi}_{L^p([-T,T];L^q)}^2\\
&\le C_T\sum_{1\le j\le j_h}(\norm{\chi_j\Op_h(a_h)\varphi}_{L^2}^2+h\norm{\varphi}_{L^2}^2)\\
&+C\sum_{1\le j\le j_h}\norm{\wtilde{\chi}_j\<x\>^{-1/2-\ep}h^{-1/2-\ep}\Op_h(a_h+b_h)e^{-itH}\varphi}_{L^2([-T,T];L^2)}^2\\
&\le C_T\norm{\varphi}_{L^2}^2+C\norm{\<x\>^{-1/2-\ep}h^{-1/2-\ep}\Op_h(a_h+b_h)e^{-itH}\varphi}_{L^2([-T,T];L^2)}^2\\
&\le C_{T,\ep}h^{-2\ep}\norm{\varphi}_{L^2}^2,
\end{align*}
which implies
\begin{align*}
\sum_{h}\norm{\Op_h(a_h)e^{-itH}f(h^2H)\varphi}_{L^p([-T,T];L^q)}^2 
&\le C_{T,\ep}\sum_{h} h^{-2\ep}\norm{f(h^2H)\varphi}_{L^2}^2\\
&\le C_{T,\ep}\norm{\<H\>^{\ep/2}\varphi}_{L^2}^2.
\end{align*}
We complete the proof.
\end{proof}

%%%%%%%%%% Bibliography %%%%%%%%%%%%%%%%%%%

\end{document}